\documentclass[aos,MSNbibl,seceqn,nameyear,dvips]{arximspdf}
\usepackage{mathbh,multirow}
\usepackage{graphicx}

\doi{10.1214/12-AOS989} 
\volume{40}
\issue{2}
\pubyear{2012}
\firstpage{1024}
\lastpage{1060}

\makeatletter
\newproclaim{remark}{Remark}[section]
\newtheorem{proposition}{Proposition}[section]
\newtheorem{theorem}{Theorem}[section]
\newtheorem{corollary}{Corollary}[section]

\newcommand{\eqref}[1]{(\ref{#1})}
\newcommand{\argmin}{\mathop{\mathrm{arg min}}}
\newcommand{\E}{{\mathbb{E}}}
\newcommand{\supp}{{\operatorname{Supp}}}
\newcommand{\im}{{\operatorname{Im}}}
\newcommand{\re}{{\operatorname{Re}}}
\newcommand{\tr}{{\operatorname{Tr}}}
\newcommand{\diag}{{\operatorname{Diag}}}
\newcommand{\bone}{\mathbh{1}}
\newcommand{\C}{{\mathbb{C}}}
\newcommand{\R}{{\mathbb{R}}}
\newcommand{\wh}{\widehat}
\newcommand{\wt}{\widetilde}
\newcommand{\gs}{\gamma}
\makeatother

\begin{document}
\begin{frontmatter}

\title{Nonlinear shrinkage estimation of large-dimensional covariance matrices}
\runtitle{Nonlinear shrinkage estimation}

\begin{aug}
\author[a]{\fnms{Olivier} \snm{Ledoit}\ead[label=e1]{olivier.ledoit@econ.uzh.ch}}
\and
\author[a]{\fnms{Michael} \snm{Wolf}\corref{}\thanksref{t1}\ead[label=e2]{michael.wolf@econ.uzh.ch}}
\runauthor{O. Ledoit and M. Wolf}
\affiliation{University of Zurich}
\address[a]{Department of Economics\\
University of Zurich \\
CH-8032 Zurich \\
Switzerland \\
\printead{e1}\\
\hphantom{E-mail: }\printead*{e2}} 
\end{aug}
\thankstext{t1}{Supported by the NCCR Finrisk project ``New Methods in
Theoretical and Empirical Asset Pricing.''}

\received{\smonth{10} \syear{2010}}
\revised{\smonth{12} \syear{2011}}

%
\begin{abstract}
Many statistical applications require an estimate of a covariance
matrix and/or
its inverse. When the matrix dimension is large compared to the sample
size, which happens frequently, the sample covariance matrix is known
to perform poorly and may suffer from ill-conditioning. There already
exists an extensive literature concerning improved estimators in such
situations. In the absence of further knowledge about
the structure of the true covariance matrix,
the most successful approach so far, arguably, has
been shrinkage estimation. Shrinking the sample covariance matrix to a
multiple of the identity, by taking a weighted average of the two,
turns out to be equivalent to linearly shrinking the sample eigenvalues to
their grand mean, while retaining the sample eigenvectors. Our paper
extends this approach by considering nonlinear transformations of the
sample eigenvalues. We show how to construct an estimator that is
asymptotically equivalent to an oracle
estimator suggested in previous work.
As demonstrated in extensive Monte Carlo simulations,
the resulting \textit{bona fide} estimator can result in sizeable
improvements over
the sample covariance matrix and also over linear shrinkage.
\end{abstract}

%
\begin{keyword}[class=AMS]
\kwd[Primary ]{62H12}
\kwd[; secondary ]{62G20}
\kwd{15A52}.
\end{keyword}
\begin{keyword}
\kwd{Large-dimensional asymptotics}
\kwd{nonlinear shrinkage}
\kwd{rotation equivariance}.
\end{keyword}

\end{frontmatter}


\section{Introduction}

Many statistical applications require an estimate of a~covariance
matrix and/or of its inverse when the matrix dimension, $p$, is large
compared to the sample size, $n$. It is well known that in such
situations, the usual estimator---the sample covariance matrix---performs
 poorly. It tends to be far from the population covariance
matrix and ill-conditioned. The goal then becomes to find estimators
that outperform the sample covariance matrix, both in finite samples
and asymptotically. For the purposes of asymptotic analyses, to
reflect the fact that $p$ is large compared to $n$, one has to employ
large-dimensional asymptotics where $p$ is allowed to go to infinity
together with $n$. In contrast, standard asymptotics would assume that
$p$ remains fixed while $n$ tends to infinity.

One way to come up with improved estimators is to incorporate
additional knowledge in the estimation process, such as sparseness, a
graph model or a factor model; for example, see
\citet{bickel:levina:2008}, \citet{rohde:tsybakov:2010}, \citet{cai:zhou:2012},
\citet{ravikumar:et:al:2008}, \citet{bala:et:al:2008},
\citet{khare:bala:2011}, \citet{fan:fan:lv:2008} and the references
therein.

However, not always is such additional knowledge available or
trustworthy. In this general case, it is reasonable to require that
covariance matrix estimators be rotation-equivariant. This means that
rotating the data by some orthogonal matrix rotates the estimator in
exactly the same way. In terms of the well-known decomposition of a
matrix into eigenvectors and eigenvalues, an estimator is
rotation-equivariant if and only if it has the same eigenvectors as
the sample covariance matrix. Therefore, it can only differentiate
itself by its eigenvalues.

\citet{ledoit:wolf:2004a} demonstrate that the largest sample
eigenvalues are systematically biased upwards, and the smallest ones
downwards. It is advantageous to correct this bias by pulling down the
largest eigenvalues and pushing up the smallest ones, toward the
grand mean of all sample eigenvalues. This is an application of the
general shrinkage principle, going back to \citet{stein:1956}. Working
under large-dimensional asymptotics, \citet{ledoit:wolf:2004a} derive
the optimal \emph{linear} shrinkage formula (when the loss is defined
as the Frobenius norm of the difference between the estimator and the
true covariance matrix). The same shrinkage intensity is applied to
all sample eigenvalues, regardless of their positions. For example, if
the linear shrinkage intensity is 0.5, then every sample eigenvalue is
moved half-way toward the grand mean of all sample eigenvalues. \citet
{ledoit:wolf:2004a} both derive asymptotic optimality properties of the
resulting estimator of the covariance matrix and demonstrate that it
has desirable finite-sample properties via simulation studies.

A cursory glance at the \citet{marcenko:pastur:1967} equation, which
governs the relationship between sample and population eigenvalues
under large-dimensional asymptotics, shows that linear shrinkage is
the first-order approximation to a fundamentally nonlinear
problem. How good is this approximation? \citet{ledoit:wolf:2004a} are
very clear about this. Depending on the situation at hand, the
improvement over the sample covariance matrix can either be gigantic
or minuscule. When $p/n$ is large, and/or the population eigenvalues
are close to one another, linear shrinkage captures most of the
potential improvement over the sample covariance matrix. In the
opposite case, that is, when $p/n$ is small and/or the population
eigenvalues are dispersed, linear shrinkage hardly improves at all
over the sample covariance matrix.

The intuition behind the present paper is that the first-order
approximation does not deliver a sufficient improvement when
higher-order effects are too pronounced. The cure is to upgrade to
\emph{nonlinear} shrinkage estimation of the covariance matrix. We get
away from the one-size-fits-all approach by applying an
individualized shrinkage intensity to every sample eigenvalue. This is
more challenging mathematically than linear shrinkage because many
more parameters need to be estimated, but it is worth the extra
effort. Such an estimator has the potential to asymptotically at least
match the linear
shrinkage estimator of \citet{ledoit:wolf:2004a} and often do a lot better,
especially when linear shrinkage does not deliver a sufficient
improvement over the sample
covariance matrix. As will be shown later in the paper, this is indeed
what we achieve here. By providing substantial improvement over the
sample covariance matrix throughout the entire parameter space,
instead of just part of it, the nonlinear shrinkage estimator is as
much of a~step forward relative to linear shrinkage as linear
shrinkage was relative to the sample covariance matrix.
In terms of finite-sample performance, the linear shrinkage estimator rarely
performs better than the nonlinear shrinkage estimator. This happens only
when the linear shrinkage estimator is (nearly) optimal
already. However, as we show in simulations, the outperformance over
the nonlinear shrinkage estimator is very small in such
cases. Most of the time, the linear shrinkage estimator is far from
optimal, and nonlinear shrinkage then offers a considerable amount of
finite-sample improvement.

A formula for nonlinear shrinkage intensities has recently been
proposed by \citet{ledoit:peche:2011}. It is motivated by a
large-dimensional asymptotic approximation to the optimal
finite-sample rotation-equivariant shrinkage formula under the
Frobenius norm. The advantage of the
formula of \citet{ledoit:peche:2011} is that it does not depend on the
unobservable population
covariance matrix: it only depends on the distribution of sample
eigenvalues. The disadvantage is that the resulting covariance matrix
estimator is an \emph{oracle} estimator in that it depends on the
``limiting'' distribution of sample eigenvalues, not the observed
one. These two objects are very different. Most critically, the
limiting empirical cumulative distribution function (c.d.f.) of sample
eigenvalues is continuously
differentiable, whereas the observed one is,
by construction, a step function.

The main contribution of the present paper is to obtain a \textit{bona
fide} estimator of the covariance matrix that is asymptotically as
good as the oracle estimator. This is done by consistently estimating
the oracle nonlinear shrinkage intensities of \citet
{ledoit:peche:2011}, in a uniform sense. As a
by-product, we also derive a new estimator of the limiting empirical
c.d.f. of population eigenvalues. A previous such estimator was
proposed by \citet{karoui:2008}.

Extensive Monte Carlo simulations indicate that our covariance matrix
estimator improves substantially over the sample covariance matrix,
even for matrix dimensions as low as $p=30$. As expected, in some
situations the nonlinear\vadjust{\goodbreak} shrinkage estimator performs as well as Ledoit
and Wolf's (\citeyear{ledoit:wolf:2004a}) linear shrinkage estimator, while in others, where
higher-order effects are more pronounced, it does substantially better.
Since the magnitude of higher-order effects depends on the population
covariance matrix, which is unobservable, it is always safer \textit{a priori}
to use nonlinear shrinkage.

Many statistical applications require an estimate of the precision matrix,
which is the inverse of the covariance matrix, instead of (or in
addition to) an estimate of the covariance matrix itself. Of course,
one possibility is to simply take the inverse of the nonlinear
shrinkage estimate of the covariance matrix itself. However, this
would be \emph{ad hoc}. The superior approach is to estimate the
inverse covariance matrix directly by nonlinearly shrinking the
inverses of the sample eigenvalues. This gives quite different and
markedly better results. We provide a detailed, in-depth solution for
this important problem as well.

The remainder of the paper is organized as
follows. Section~\ref{sec:asym} defines our framework for
large-dimensional asymptotics and reviews some fundamental results
from the corresponding literature. Section~\ref{sec:oracle} presents
the oracle shrinkage estimator that motivates our \textit{bona fide}
nonlinear shrinkage estimator. Sections~\ref{sec:m} and~\ref{sec:H}
show that the \textit{bona fide} estimator is consistent for the oracle
estimator. Section~\ref{sec:MC} examines finite-sample
behavior via Monte Carlo simulations. Finally, Section~\ref
{sec:conclusions} concludes. All
mathematical proofs are collected in the supplement [\citet
{ledoit:wolf:supplement}].

\section{Large-dimensional asymptotics}\label{sec:asym}

\subsection{Basic framework}
Let $n$ denote the sample size and $p \equiv p(n)$ the number of variables,
with $p/n\to c\in(0,1)$ as $n\to\infty$.
This framework is known as
large-dimensional asymptotics.
The restriction to the case $c < 1$ that we make here somewhat
simplifies certain mathematical results as well as the implementation
of our routines in software. The case $c > 1$, where the sample
covariance matrix is singular, could be handled by
similar methods, but is left to future research.

The following set of assumptions will be maintained throughout the
paper.

\begin{enumerate}[(A3)]
\item[(A1)] The population covariance matrix $\Sigma_n$ is a nonrandom
$p$-dimensional positive definite matrix.

\item[(A2)] Let $X_n$ be an $n\times p$ matrix of
real independent and identically distributed (i.i.d.) random
variables with zero mean and unit variance. One only observes $Y_n
\equiv
X_n \Sigma_n^{1/2}$, so neither $X_n$ nor $\Sigma_n$ are
observed on their own.

\item[(A3)] Let
$((\tau_{n,1},\ldots,\tau_{n,p});(v_{n,1},\ldots,v_{n,p}))$ denote
a system of
eigenvalues and eigenvectors of $\Sigma_n$. The
empirical distribution function (e.d.f.) of the population
eigenvalues is defined as $\forall t\in\R,
H_n(t)\equiv p^{-1}\sum_{i=1}^p\bone_{[\tau_{n,i},+\infty)}(t)$,
where $\bone$
denotes the indicator function of a set. We assume $H_n(t)$
converges to some limit $H(t)$ at all points of continuity of $H$.

\item[(A4)] $\supp(H)$, the support of $H$, is the union of a
finite number of closed intervals, bounded away from zero and
infinity. Furthermore, there exists a compact interval in $(0, +
\infty)$ that contains $\supp(H_n)$ for all $n$ large enough.
\end{enumerate}

Let $((\lambda_{n,1},\ldots,\lambda_{n,p});(u_{n,1},\ldots
,u_{n,p}))$ denote a
system of eigenvalues and eigenvectors of the sample covariance matrix
$S_n \equiv n^{-1} Y_n^\prime Y_n =
n^{-1}\Sigma_n^{1/2}X_n'X_n\*\Sigma_n^{1/2}$. We can assume that the
eigenvalues are sorted in increasing order without loss of
generality (w.l.o.g.).
The first subscript, $n$, will be omitted when no confusion is
possible. The e.d.f. of the sample eigenvalues is defined as
$\forall\lambda\in\R,
F_n(\lambda)\equiv
p^{-1}\sum_{i=1}^p\bone_{[\lambda_i,+\infty)}(\lambda)$.

In the remainder of the paper, we shall use the notation $\re(z)$ and
$\im(z)$ for the real and imaginary parts, respectively, of a complex
number $z$, so
that
\[
\forall z \in\mathbb{C}\qquad z = \re(z) + i \cdot\im(z).
\]

The Stieltjes transform of a nondecreasing function $G$ is defined by
%
\begin{equation} \label{e:stieltjes}
\forall z\in\mathbb{C}^+\qquad m_G(z) \equiv
\int_{-\infty}^{+\infty}\frac{1}{\lambda-z}\,dG(\lambda),
\end{equation}
where $\C^+$ is the half-plane of complex numbers with strictly
positive imaginary part. The Stieltjes transform has a well-known
inversion formula,\vspace*{-1pt}
\[
G(b)-G(a)=\lim_{\eta\to0^+}\frac{1}{\pi}\int_a^b\im
[m_G(\xi+i\eta) ]\,d\xi,\vspace*{-1pt}
\]
which holds if $G$ is continuous at $a$ and $b$. Thus, the Stieltjes
transform of the e.d.f. of sample eigenvalues is\vspace*{-1pt}
\[
\forall z\in\mathbb{C}^+\qquad
m_{F_n}(z)=\frac{1}{p}\sum_{i=1}^p\frac{1}{\lambda_i-z}
=\frac{1}{p}\tr[(S_n-zI)^{-1} ],\vspace*{-1pt}
\]
where $I$ denotes a conformable identity matrix.\vspace*{-1pt}

\subsection{Mar\protect{\v{c}}enko--Pastur equation and reformulations}
\label{ss:mp}

Mar{\v{c}}enko and
Pastur (\citeyear{marcenko:pastur:1967}) and others have proven that $F_n(\lambda)$
converges almost surely (a.s.) to some nonrandom limit $F(\lambda)$ at
all points of
continuity of $F$ under certain sets of assumptions. Furthermore,
Mar\v{c}enko and Pastur discovered the equation that relates $m_F$ to
$H$. The most convenient expression of the Mar\v{c}enko--Pastur
equation is the one found in
Silverstein [(\citeyear{silverstein:1995}), equation (1.4)],
%
\begin{equation}
\label{eq:MP}
\forall z\in\mathbb{C}^+\qquad
m_F(z)=\int_{-\infty}^{+\infty}\frac{1}{\tau
[1-c-c z m_F(z) ]-z}\, dH(\tau).
\end{equation}
This version of the Mar\v{c}enko--Pastur equation is the one that we
start out with. In addition, \citet{silverstein:choi:1995} showed
that
\[
\forall\lambda\in\R-\{0\}\qquad \lim_{z\in\mathbb{C}^+\to\lambda
} m_F(z)\equiv
\breve{m}_F(\lambda)\vadjust{\goodbreak}
\]
 exists, and that $F$ has a continuous derivative
$F'=\pi^{-1}\im[\breve{m}_F ]$ on all of $\R$ with $F^\prime
\equiv0$ on $(-\infty, 0]$. For purposes that will become apparent
later, it is useful to reformulate the Mar\v{c}enko--Pastur equation.

The limiting e.d.f. of the eigenvalues of $n^{-1} Y_n'Y_n =
n^{-1} \Sigma_n^{1/2} X_n^\prime X_n \Sigma_n^{1/2}$ was defined as
$F$. In addition, define the
limiting e.d.f. of the eigenvalues of
$n^{-1}Y_n Y_n' = n^{-1} X_n\Sigma_nX_n'$ as
$\underline{F}$. It then holds
\begin{eqnarray*}
\forall x\in\R\qquad \underline{F}(x)&=&(1-c) \bone_{[0,+\infty)}(x)+c
F(x),
\\
\forall x\in\R\qquad F(x)&=&\frac{c-1}{c}\bone_{[0,+\infty)}(x)+\frac
{1}{c} \underline{F}(x),
\\
\forall z \in\C^+\qquad m_{\underline{F}}(z)&=&\frac{c-1}{z}+c
m_F(z),
\\
\forall z \in\C^+\qquad m_F(z)&=&\frac{1-c}{c z}+\frac{1}{c} m_{\underline
{F}}(z).
\end{eqnarray*}
With this notation, equation (1.3) of \citet{silverstein:choi:1995}
rewrites the Mar\v{c}enko--Pastur equation in the following way: for
each $z\in\C^+$, $m_{\underline{F}}(z)$ is the unique solution in
$\C^+$ to the equation
%
\begin{equation}
\label{eq:mp}
m_{\underline{F}}(z)=- \biggl[
{z-c\int_{-\infty}^{+\infty}\frac{\tau}{1+\tau
m_{\underline{F}}(z)}\,dH(\tau)} \biggr]^{-1}.
\end{equation}
Now introduce $u_{\underline{F}}(z) \equiv-1/{m_{\underline{F}}(z)}$.
Notice that $u_{\underline{F}}(z)\in\C^+\iff m_{\underline
{F}}(z)\in\C^+$. The mapping from $u_{\underline{F}}(z)$ to
$m_{\underline{F}}(z)$ is one-to-one on $\C^+$.

With this change of variable, equation (\ref{eq:mp}) is equivalent to
saying that for each $z\in\C^+$, $u_{\underline{F}}(z)$ is the
unique solution in $\C^+$ to the equation
%
\begin{equation}
\label{eq:mp2}
u_{\underline{F}}(z)=z+c u_{\underline{F}}(z)\int_{-\infty
}^{+\infty}\frac{\tau}{\tau-u_{\underline{F}}(z)}\,dH(\tau).
\end{equation}
Let the linear operator $L$ transform any c.d.f. $G$ into
\[
LG(x) \equiv\int_{-\infty}^{x}\tau \,dG(\tau).
\]
Combining $L$ with the Stieltjes transform, we get
\[
m_{LG}(z)=\int_{-\infty}^{+\infty}\frac{\tau}{\tau-z}\,dG(\tau)=1+z
m_G(z).
\]
Thus, we can rewrite equation (\ref{eq:mp2}) more concisely as
%
\begin{equation}
\label{eq:mp3}
u_{\underline{F}}(z)=z+c u_{\underline{F}}(z)
m_{LH} (u_{\underline{F}}(z) ).
\end{equation}
As Silverstein and Choi [(\citeyear{silverstein:choi:1995}), equation (1.4)] explain, the function
defined in equation (\ref{eq:mp}) is invertible. Thus we can define
the inverse function
%
\begin{equation}
\label{eq:invertm}
z_{\underline{F}}(m) \equiv-\frac{1}{m}+c\int_{-\infty}^{+\infty
}\frac{\tau}{1+\tau
m}\,dH(\tau).
\end{equation}
We can do the same thing for equation (\ref{eq:mp3}) and define the
inverse function
%
\begin{equation}
\label{eq:invertu}
\widetilde{z}_{\underline{F}}(u) \equiv u-c u m_{LH}(u).
\end{equation}
Equations \eqref{eq:MP},
(\ref{eq:mp}), (\ref{eq:mp3}), (\ref{eq:invertm}) and
(\ref{eq:invertu}) are all completely equivalent to one another;
solving any one of them means having solved them all. They are all
just reformulations of the Mar\v{c}enko--Pastur equation.

As will be detailed in Section~\ref{sec:oracle}, the oracle nonlinear
shrinkage estimator
of $\Sigma_n$ involves the quantity ${\breve m}_F(\lambda)$, for
various inputs
$\lambda$. Section~\ref{ss:solving} describes how this quantity
can be
found in the hypothetical case that $F$ and $H$ are actually known.
This will then allow us later to discuss consistent estimation of
${\breve m}_F(\lambda)$ in the realistic case when $F$ and $H$ are unknown.

\subsection{Solving the Mar\protect{\v{c}}enko--Pastur equation}
\label{ss:solving}

Silverstein and Choi (\citeyear{silverstein:choi:1995})
 explain how the support of $F$, denoted by $\supp(F)$, is
determined. Let $B \equiv\{u \in\R\dvtx u \neq0, u \in
\supp^\complement(H) \}$. Then plot the function
$\tilde z_{\underline F}(u)$ of \eqref{eq:invertu}
on the set $B$. Find the extreme values on each interval. Delete these
points and everything in between on the real line. Do this for all
increasing intervals. What is left is just $\supp(F)$; see Figure~1 of
\citet{bai:silverstein:1998} for an illustration.

To simplify, we will assume from here on that
$\supp(F)$ is a single compact interval, bounded away from zero,
with $F^\prime> 0$ in the interior of this interval.
But if $\supp(F)$ is the union of a finite
number of such intervals, the arguments presented in this section as
well as in the remainder of the paper apply separately to each
interval. In particular, our consistency results presented in
subsequent sections can be easily extended to this more general
case. On the other hand, the even
more general case of $\supp(F)$ being the union of an infinite number
of such intervals or being a noncompact interval is ruled out by
assumption (A4).
By our assumption then, $\supp(F)$ is given by the compact interval
$[\widetilde{z}_{\underline{F}}(u_1),\widetilde{z}_{\underline{F}}(u_2)]$
for some $u_1 < u_2$. To
keep the notation shorter in what follows, let $\wt z_1 \equiv\wt
z_{\underline F}(u_1)$ and
$\wt z_2 \equiv\wt z_{\underline F}(u_2)$.

We know that for every
$\lambda$ in the interior of $\supp(F)$,
there exists a unique $v\in\C^+$, denoted by $v_\lambda$, such that
%
\begin{equation}\label{e:cum}
v_\lambda-c v_\lambda m_{LH}(v_\lambda)=\lambda.
\end{equation}
We further know that
\[
F'(\lambda)=\frac{1}{c}\underline{F}'(\lambda)=\frac{1}{c\pi}\im
[\breve
m_{\underline{F}}(\lambda)]
=\frac{1}{c\pi}\im\biggl[-\frac{1}{v_\lambda} \biggr].
\]

The converse is also true.
Since $\supp(F) =
[\widetilde{z}_{\underline{F}}(u_1),\widetilde{z}_{\underline{F}}(u_2)]$,
for every $x\in(u_1,u_2)$, there exists a
unique $y > 0$, denoted by $y_x$, such that
\[
(x+iy_x)-c (x+iy_x) m_{LH}(x+iy_x)\in\R.
\]
In other words, $y_x$ is the unique value of $y > 0$ for which
$\im[(x+iy)-c(x+iy)m_{LH}(x+iy) ]=0$. Also, if $\lambda_x$ denotes
the value of
$\lambda$ for which we have $(x+iy_x)-c (x+iy_x) m_{LH}(x+iy_x) =
\lambda$, then, by definition, $z_{\lambda_x} = x + i y_x$.

Once we find a way to
consistently estimate $y_x$ for any $x\in[u_1,u_2]$, then we have
an estimate of the (asymptotic) solution to the Mar\v{c}enko--Pastur equation.
For example,
$\im[-1/(x+iy_x) ]/(c\pi)$ is the value of the density~$F^\prime$
evaluated at $\re[(x+iy_x)-c (x+iy_x) m_{LH}(x+iy_x) ]
=(x+iy_x)-c (x+iy_x) m_{LH}(x+iy_x)$.

From the above arguments, it follows that
%
\begin{equation}\label{e:m-breve}\qquad
\forall\lambda\in(\wt z_1, \wt z_2)\qquad
\breve m_{\underline F} (\lambda) = -\frac{1}{v_\lambda}
\quad\mbox{and so}\quad
\breve m_{F} (\lambda) = \frac{1 - c}{c \lambda} - \frac{1}{c}
\frac{1}{v_\lambda}.
\end{equation}

\section{Oracle estimator}\label{sec:oracle}

\subsection{Covariance matrix}\label{ss:cov-mat}

In the absence of specific information about the true covariance
matrix $\Sigma_n$, it appears reasonable to restrict attention to the
class of estimators that are equivariant with respect to rotations of
the observed data. To be more specific, let $W$ be an arbitrary
$p$-dimensional orthogonal matrix. Let $\widehat\Sigma_n \equiv
\widehat\Sigma_n(Y_n)$
be an estimator of $\Sigma_n$. Then the
estimator is said to be \textit{rotation-equivariant} if
it satisfies $\widehat\Sigma_n(Y_n W) = W^\prime\widehat\Sigma
_n(Y_n) W$.
In other words, the estimate based on the rotated data equals the
rotation of the estimate based on the original data.
The class of rotation-equivariant estimators of the covariance matrix
is constituted of all the estimators that have the same eigenvectors
as the sample covariance matrix; for example, see
\citeauthor{perlman:2007} [(\citeyear{perlman:2007}), Section~5.4]. Every
rotation-equivariant estimator is thus of
the form
\[
U_n D_n U_n^\prime\qquad\mbox{where }
D_n \equiv\diag(d_1, \ldots, d_p) \mbox{ is diagonal},
\]
and where $U_n$ is the matrix whose $i$th column is the sample
eigenvector $u_i \equiv u_{n,i}$. This is the class we consider.

The starting objective is to find the matrix in this class that is
closest to~$\Sigma_n$. To measure distance, we choose the
Frobenius norm defined as
%
\begin{equation} \label{e:fro}\qquad
\Vert A\Vert \equiv\sqrt{\tr(A
A^\prime)/r}\qquad \mbox{for any matrix $A$ of dimension $r \times m$}.
\end{equation}
[Dividing by the dimension of the square matrix $A A^\prime$ inside
the root is not
standard, but we do this for asymptotic purposes so that the Frobenius
norm remains constant equal to one for the identity matrix regardless
of the
dimension; see \citet{ledoit:wolf:2004a}.]
As a result, we end up with the following
minimization problem:
\[
\min_{D_n} \Vert U_n D_n U_n^\prime- \Sigma_n\Vert.
\]
Elementary matrix algebra shows that its solution is
%
\begin{equation}\label{e:star}\qquad
D_n^* \equiv\diag(d_1^*, \ldots, d_p^*)
\qquad\mbox{where }
d_i^* \equiv u_i^\prime\Sigma_n u_i \mbox{ for } i = 1, \ldots, p.\vadjust{\goodbreak}
\end{equation}
The interpretation of $d_i^*$ is that it captures how the
$i$th sample eigenvector~$u_i$ relates to the population covariance
matrix $\Sigma_n$ as a whole.
As a result, the finite-sample optimal estimator is given by
%
\begin{equation}\label{e:optimal-fs}
S_n^* \equiv U_n D_n^* U_n^\prime\qquad\mbox{where }
D_n^* \mbox{ is defined as in \eqref{e:star}}.
\end{equation}

By generalizing the Mar\v{c}enko--Pastur equation (\ref{eq:MP}),
\citet{ledoit:peche:2011} show that $d_i^*$ can be approximated by the
quantity
%
\begin{equation}\label{e:oracle0}
d_i^{or} \equiv\frac{\lambda_i}{ |1 - c -c \lambda_i \breve
m_F(\lambda_i) |^2} \qquad\mbox{for } i = 1, \ldots, p,
\end{equation}
from which they deduce their oracle estimator
%
\begin{equation}\label{e:oracle}
S_n^{or} \equiv U_n D_n^{or} U_n^\prime
\qquad\mbox{where }
D_n^{or} \equiv\diag(d_1^{or}, \ldots, d_p^{or}).
\end{equation}
The key difference between $D^*_n$ and $D^{or}_n$ is that the former
depends on the unobservable population covariance matrix, whereas the
latter depends on the limiting distribution of sample eigenvalues,
which makes it amenable to estimation, as explained below.

Note that $S_n^{or}$ constitutes a nonlinear shrinkage estimator: since
the value of the denominator of $d_i^{or}$ varies with $\lambda_i$, the
shrunken eigenvalues $d_i^{or}$ are obtained by applying a nonlinear
transformation to the sample eigenvalues~$\lambda_i$; see
Figure~\ref{fig:nonlinearity} for an illustration.
\citet{ledoit:peche:2011}
also illustrate in some (limited) simulations that this oracle
estimator can provide a magnitude of improvement over the
linear shrinkage estimator of \citet{ledoit:wolf:2004a}.

\subsection{Precision matrix}\label{ss:oracle:precision}

Often times an estimator of the inverse of the covariance matrix, or
the precision matrix, $\Sigma_n^{-1}$ is required. A reasonable
strategy would be to first estimate $\Sigma_n$, and to then simply take
the inverse of the resulting estimator. However, such a strategy will
generally not be optimal.

By arguments analogous to those leading up to \eqref{e:optimal-fs},
among the class of rotation-equivariant estimators, the finite-sample
optimal estimator of $\Sigma_n^{-1}$ with respect to the Frobenius
norm is given by
%
\begin{equation}\label{e:optimal-fs-precision}
P_n^* \equiv U_n A_n^* U_n^\prime\qquad\mbox{where }
a_i^* \equiv u_i^\prime\Sigma_n^{-1} u_i \mbox{ for } i = 1, \ldots
, p.
\end{equation}
In particular, note that $P_n^* \neq(S_n^*)^{-1}$ in general.

Studying the asymptotic behavior of the diagonal matrix $A_n^*$ led
\citet{ledoit:peche:2011} to the following oracle estimator:
%
\begin{eqnarray}\label{e:oracle-precision}
P_n^{or} \equiv U_n A_n^{or} U_n^\prime\nonumber
\\[-8pt]
\\[-8pt]
\eqntext{\mbox{where }
a_i^{or} \equiv\lambda_i^{-1} \bigl(1 - c - 2 c \lambda_i
\re[\breve m_F(\lambda_i)] \bigr)
\mbox{ for } i = 1, \ldots, p.}
\end{eqnarray}
In particular, note that $P_n^{or} \neq(S_n^{or})^{-1}$ in general.

\begin{remark}
One can see that both oracle estimators $S_n^{or}$ and
$P_n^{or}$ involve the unknown quantities $\breve
m_F(\lambda_i)$, for $i = 1, \ldots, p$. As a result, they are not
\textit{bona fide} estimators.
However, being able to consistently estimate~$\breve
m_F(\lambda)$, uniformly in~$\lambda$, will allow us to construct
\textit{bona fide} estimators $\widehat S_n$ and $\widehat P_n$
that converge to their respective oracle counterparts
almost surely (in the sense that the Frobenius norm of the difference
converges to zero almost surely).

Section~\ref{sec:m} explains how to construct a uniformly consistent
estimator of~$\breve m_F(\lambda)$
based on a consistent estimator of
$H$, the limiting spectral distribution of the population eigenvalues.
Section~\ref{sec:H} discusses how to construct a~consistent
estimator of~$H$ from the data.
\end{remark}

\subsection{\texorpdfstring{Further details on the results of Ledoit and P\'{e}ch\'{e} (\citeyear{ledoit:peche:2011})}{Further details on the results of Ledoit and P\'{e}ch\'{e} (2011)}}

Ledoit and P\'{e}ch\'{e} (\citeyear{ledoit:peche:2011}) (hereafter LP) study functionals of the type
%
\begin{eqnarray}\label{eq:general}
\forall z\in\mathbb{C}^+\qquad
\Theta^g_N(z)&\equiv&\frac{1}{N}\sum_{i=1}^N\frac{1}{\lambda
_i-z}\sum_{j=1}^N |u_i^*v_j |^2\times g(\tau_j)\nonumber
\\[-8pt]
\\[-8pt]
&=&\frac{1}{N}\tr[(S_N-zI)^{-1}g(\Sigma_N) ],\nonumber
\end{eqnarray}
%
where $g$ is any real-valued univariate function satisfying suitable
regularity conditions. Comparison with equation \eqref{e:stieltjes}
reveals that this family of functionals generalizes the Stieltjes
transform, with the Stieltjes transform corresponding to the special
case $g\equiv1$. What is of interest is what happens for other,
nonconstant functions $g$.

It turns out that it is possible to generalize the Mar\v{c}enko--Pastur
result~(\ref{eq:MP}) to any function $g$ with finitely many points of
discontinuity. Under assumptions that are usual in the Random Matrix
Theory literature, LP prove in their Theorem~2 that there exists a
nonrandom function
$\Theta^g$ defined over $\mathbb{C}^+$ such that $\Theta_N^g(z)$
converges a.s. to $\Theta^g(z)$ for all $z\in\mathbb{C}^+$.
Furthermore,~$\Theta^g$ is given by
%
\begin{equation}
\label{eq:theta}
\forall z\in\mathbb{C}^+\qquad
\Theta^g(z) \equiv\int_{-\infty}^{+\infty}\frac{g(\tau)}{\tau
[1-c-c z m_F(z) ]-z}\, dH(\tau).
\end{equation}
What is remarkable is that, as one moves from the constant
function $g\equiv1$ to any other function $g(\tau)$, the integration
kernel $\frac{g(\tau)}{\tau
[1-c-c z m_F(z) ]-z}$ remains\vspace*{1pt} unchanged. Therefore
equation (\ref{eq:theta}) is a direct generalization of Mar\v{c}enko
and Pastur's foundational result.

The power and usefulness of this generalization become apparent once
one starts plugging specific, judiciously chosen functions $g(\tau)$
into equation (\ref{eq:theta}). For the purpose of illustration, LP
work out three examples of functions $g(\tau)$.

The first example of LP is
$g(\tau) \equiv\bone_{(-\infty,\tau)}$, where $\bone$ denotes the
indicator function of a set.\vadjust{\goodbreak} It enables them to characterize the
asymptotic location of sample eigenvectors relative to population
eigenvectors. Since this result is not directly relevant to the
present paper, we will not elaborate further, and refer the interested
reader to LP's Section~1.2.

The second example of LP is $g(\tau)
\equiv\tau$. It
enables them to characterize the asymptotic behavior of the quantities
$d_i^{or}$ introduced in equation (\ref{e:oracle0}). More formally,
for any $u\in(0,1)$, define
%
\begin{equation}
\Delta_{n}^*(u)\equiv\frac{1}{p}\sum_{i=1}^{ \lfloor u\cdot p
\rfloor}d_i^* \quad\mbox{and}\quad
\Delta_{n}^{or}(u)\equiv\frac{1}{p}\sum_{i=1}^{ \lfloor u\cdot p
\rfloor}d_i^{or},
\end{equation}
where $\lfloor\cdot\rfloor$ denotes the integer part.
LP's Theorem~4 proves that
$\Delta_{n}^*(u)-\Delta_{n}^{or}(u) \rightarrow0$ a.s.

The third example of LP is $g(\tau) \equiv1/\tau$. It enables them
to characterize the asymptotic behavior of the quantities $a_i^{or}$
introduced in equation (\ref{e:oracle-precision}). For any $u\in(0,1)$ define
%
\begin{equation}
\Psi_{n}^*(u)\equiv\frac{1}{p}\sum_{i=1}^{ \lfloor u\cdot p \rfloor
}a_i^* \quad\mbox{and}\quad
\Psi_{n}^{or}(u)\equiv\frac{1}{p}\sum_{i=1}^{ \lfloor u\cdot p
\rfloor}a_i^{or}.
\end{equation}
LP's Theorem~5 proves that $\Psi_{n}^*(u)-\Psi_{n}^{or}(u)
\rightarrow0$ a.s.

\section{\texorpdfstring{Estimation of $\breve m_F(\lambda)$}{Estimation of mF(lambda)}}\label{sec:m}

Fix $x \in[u_1 + \eta, u_2 - \eta]$, where $\eta> 0$ is some small number.
From the previous discussion in Section~\ref{sec:asym}, it follows that
the equation
\[
\im[x + i y - c (x + iy) m_{L H}(x + iy) ] = 0
\]
has a unique solution $y \in(0, +\infty)$, called $y_x$.
Since $u_1 < x < u_2$, it follows that $y_x > 0$; for $x =
u_1$ or $x = u_2$, we would have $y_x = 0$ instead.
The goal is to consistently estimate $y_x$,
uniformly in $x \in[u_1 + \eta, u_2 - \eta]$.

Define for any c.d.f. $G$ and for any $d > 0$, the real function
\[
g_{G,d}(y, x) \equiv| \im[x + iy - d (x + iy) m_{LG}(x +
iy) ] |.
\]
With this notation, $y_x$ is the unique minimizer in $(0, +\infty)$
of $g_{H,c}(y, x)$ then. In particular, $g_{H,c}(y_x, x) = 0$.

In the remainder of the paper, the symbol $\Rightarrow$ denotes weak
convergence (or
convergence in distribution).
\begin{proposition}\label{prop:y0-uni}
\emph{(i)}
Let $\{\widehat H_n\}$ be a sequence of
probability measures with $\widehat H_n \Rightarrow H$.
Let $\{\widehat c_n\}$ be a sequence of positive real numbers with
$\wh c_n \to c$. Let $K \subseteq(0, \infty)$ be a compact interval
satisfying $ \{y_x \dvtx x \in[u_1 + \eta, u_2 - \eta] \} \subseteq K$.
For a given $x \in[u_1 + \eta, u_2 - \eta]$,
let $\wh y_{n, x} \equiv\min_{y \in K}
g_{\widehat{H}_n, \wh c_n}(y, x)$.
It then holds that $\wh y_{n, x} \to y_x$ uniformly in $x
\in[u_1 + \eta, u_2 - \eta]$.

\emph{(ii)} In case of $\widehat H_n \Rightarrow H$ a.s., it
holds that $\widehat y_{n, x} \to y_x $ a.s. uniformly in $x
\in[u_1 + \eta, u_2 - \eta]$.\vadjust{\goodbreak}
\end{proposition}

It should be pointed out that the assumption $ \{y_x \dvtx x \in[u_1 +
\eta, u_2 - \eta] \} \subseteq K$ is not really restrictive, since one
can choose $K \equiv[\varepsilon, 1/\varepsilon]$, for $\varepsilon
$ arbitrarily small.

We also need to solve the ``inverse'' estimation problem, namely
starting with $\lambda$ and recovering the corresponding $v_\lambda$.
Fix $\lambda\in[\wt z_1 + \wt\delta, \wt z_2 - \wt\delta]$,
where $\wt\delta> 0$ is some small number.
From the previous discussion, it follows that
the equation
\[
v - c v m_{LH}(v) = \lambda
\]
has a unique solution $v \in\C^+$, called $v_\lambda$.
The goal is to consistently estimate~$v_{\lambda}$, uniformly in
$\lambda\in[\wt z_1 + \wt\delta, \wt z_2 - \wt\delta]$.

Define for any c.d.f. $G$ and for any $d > 0$, the real function
\[
h_{G,d}(v, \lambda) \equiv|v - d v m_{LG}(v) - \lambda
|.
\]
With this notation, $v_\lambda$ is then the unique minimizer in $\C^+$
of $h_{H,c}(v, \lambda)$. In particular, $h_{H,c}(v_\lambda, \lambda
) = 0$.

\begin{proposition}\label{prop:zl-uni}
\emph{(i)}
Let $\{\widehat H_n\}$ be a sequence of
probability measures with $\widehat H_n \Rightarrow H$.
Let $\{\widehat c_n\}$ be a sequence of positive real numbers with
$\wh c_n \to c$.
Let $K \subseteq\C^+$ be a compact set
satisfying $ \{v_\lambda\dvtx \lambda\in[\widetilde z_1 + \widetilde
\delta, \widetilde z_2 - \widetilde\delta] \} \subseteq K$.
For a given $\lambda\in[\wt z_1 + \wt\delta, \wt z_2 - \wt\delta]$,
let $\wh v_{n, \lambda} \equiv\min_{v \in K}
h_{\widehat{H}_n, \wh c_n}(v, \lambda)$.
It then holds that $\wh v_{n, \lambda} \to v_\lambda$ uniformly in
$\lambda
\in[\wt z_1 + \wt\delta, z_2 - \wt\delta]$.

\emph{(ii)} In case of $\widehat H_n \Rightarrow H$ a.s., it
holds that $\widehat v_{n, \lambda} \to v_\lambda$ a.s. uniformly in
$\lambda
\in[\wt z_1 + \wt\delta, z_2 - \wt\delta]$.
\end{proposition}

Being able to find consistent estimators of $v_\lambda$, uniformly in
$\lambda$, now
allows us to find consistent estimators of $\breve m_F(\lambda)$,
uniformly in $\lambda$, based on \eqref{e:m-breve}.
Our estimator of $\breve m_F(\lambda)$ is given by
%
\begin{equation}\label{e:m-breve-hat}
\breve m_{F_{\widehat H_n, \widehat c_n}}(\lambda) \equiv
\frac{1 - \widehat c_n}{\widehat c_n \lambda}
- \frac{1}{\widehat c_n}
\frac{1}{ \widehat v_{n, \lambda}}.
\end{equation}
This, in turn, provides us with a consistent estimator of
$S_n^{or}$, the oracle nonlinear shrinkage estimator of
$\Sigma_n$. Define
%
\begin{eqnarray}\label{e:nl-estimator}
\widehat S_n \equiv U_n \widehat D_n U_n^\prime\nonumber
\\[-8pt]
\\[-8pt]
\eqntext{\mbox{where }
\widehat d_i \equiv\displaystyle\frac{\lambda_i}{ |1 - \widehat c_n - \widehat
c_n \lambda_i
\breve m_{F_{\widehat H_n, \widehat c_n}}(\lambda_i) |^2}
\mbox{ for } i = 1, \ldots, p.}
\end{eqnarray}

It also provides us with a consistent estimator of
$P_n^{or}$, the oracle nonlinear shrinkage estimator of
$\Sigma_n^{-1}$. Define
%
\begin{eqnarray}\label{e:nl-estimator-precision}
\widehat P_n \equiv U_n \widehat A_n U_n^\prime\nonumber
\\[-8pt]
\\[-8pt]
\eqntext{\mbox{where }
\widehat a_i \equiv\lambda_i^{-1} \bigl(1 - \wh c_n - 2 \wh c_n
\lambda_i \re[\breve m_{F_{\wh H_n, \wh c_n}} (\lambda_i)] \bigr)
\mbox{ for } i = 1, \ldots, p.}
\end{eqnarray}
In particular, note that $\widehat P_n \neq\widehat S_n^{-1}$ in general.

\begin{proposition}\label{prop:breve-m-uni}
\begin{enumerate}[(ii)]
\item[(i)]
Let $\{\widehat H_n\}$ be a sequence of
probability measures with $\widehat H_n \Rightarrow H$.
Let $\{\widehat c_n\}$ be a sequence of positive real numbers with
$\wh c_n \to c$.
It then holds that:
\begin{longlist}[(b)]
\item[(a)] $\breve m_{F_{\widehat H_n, \widehat c_n}}(\lambda)\to
\breve{m}_F(\lambda)$ uniformly in $\lambda\in[\wt z_1 + \wt\delta,
\wt z_2 - \wt\delta]$;
\item[(b)] $\Vert\widehat S_n - S_n^{or}\Vert \to0$;
\item[(c)] $\Vert\widehat P_n - P_n^{or}\Vert \to0$.
\end{longlist}
\item[(ii)] In case of $\widehat H_n \Rightarrow H$ a.s., it
holds that:
\begin{longlist}[(b)]
\item[(a)] $\breve m_{F_{\widehat H_n, \widehat c_n}}(\lambda)\to
\breve{m}_F(\lambda)$ uniformly in $\lambda\in[\wt z_1 + \wt\delta,
\wt z_2 - \wt\delta]$ a.s.;
\item[(b)] $\Vert\widehat S_n - S_n^{or}\Vert \to0$ a.s.;
\item[(c)] $\Vert\widehat P_n - P_n^{or}\Vert \to0$ a.s.
\end{longlist}
\end{enumerate}
\end{proposition}

\section{Estimation of $H$}
\label{sec:H}

As described before, consistent estimation of the oracle estimators of
\citet{ledoit:peche:2011} requires (uniformly) consistent
estimation of $\breve{m}_F(\lambda)$.
Since $\im[\breve{m}_F(\lambda) ]=\pi
F'(\lambda)$, one possible approach could be to take an
off-the-shelf density estimator for $F'$, based on the observed
sample eigenvalues~$\lambda_i$. There exists a large
literature on density estimation; for example, see \citet
{silverman:1986}. The real
part of $\breve{m}_F(\lambda_i)$ could be estimated in a similar
manner.

However, the sample eigenvalues do not satisfy any of the regularity
conditions usually invoked
for the underlying data. It really is not clear at all whether an
off-the-shelf density estimator applied
to the sample eigenvalues would result in consistent estimation of $F'$.

Even if this issue was somehow resolved, using such a generic procedure
would not exploit the specific features of the
problem. Namely: $F$ is not just any distribution; it is a
distribution of sample eigenvalues. It is the solution to the
Mar\v{c}enko--Pastur equation for some $H$. This is valuable
information that narrows down considerably the set of possible
distributions~$F$. Therefore an estimation procedure specifically
designed to incorporate this a priori knowledge would be better
suited to the problem at hand. This is the approach we select.

In a nutshell: our estimator of $F$ is the c.d.f. that is closest to
$F_n$ among the c.d.f.'s that are a solution to the Mar\v{c}enko--Pastur
equation for some~$\widetilde{H}$ and for $\widetilde c \equiv
\widehat c_n \equiv p/n$.
The ``underlying'' distribution
$\widetilde{H}$ that produces the thus obtained estimator of $F$ is, in
turn, our estimator of $H$. If we can show that this estimator of
$H$ is consistent, then the results of the previous section demonstrate
that the implied estimator of $\breve m_F(\lambda)$ is uniformly
consistent.

Section~\ref{ss:H-cons} derives theoretical properties of this
approach, while Section~\ref{sec:procedure} discusses various issues
concerning the practical implementation.

\subsection{Consistency results}\label{ss:H-cons}

\hspace*{-1.365pt}For a grid of real numbers $Q \equiv\{\ldots, t_{-1}, t_0,\allowbreak t_1,
\ldots\}\subseteq\R$, with $t_{k-1} < t_k$,
define the corresponding
grid size $\gs$ as
\[
\gs\equiv\sup_k (t_{k} - t_{k-1}).
\]
A grid $Q$ is said to cover a compact interval $[a, b] \subseteq\R$ if
there exists at least one $t_k \in Q$ with $t_k \le a$ and at least
another $t_{k^\prime} \in Q$ with $b \le t_{k^\prime}$.
A~sequence of grids $\{Q_n\}$ is said to eventually cover a compact
interval $[a, b]$ if for every $\phi> 0$ there exist
$N \equiv N(\phi) $ such that $Q_n$ covers
the compact interval $[a + \phi, b-\phi]$
for all $n \ge N$.

For any probability measure $\widetilde H$ on the real line
and for any $\wt c > 0$, let $F_{\widetilde H, \wt c}$ denote
the c.d.f. on the real line induced by the corresponding solution of
the Mar\v{c}enko--Pastur
equation.
More specifically, for each $z \in\C^+$, $m_{F_{\widetilde H,\wt c}}(z)$
is the unique solution for $m \in\C^+$ to the equation
\[
m =
\int_{-\infty}^{+\infty}\frac{1}{\tau
[1-\wt c- \wt c z m ]-z} \,d \widetilde H(\tau).
\]
In this notation, we then have $F = F_{H,c}$.

It follows from \citet{silverstein:choi:1995} again
that
\[
\forall\lambda\in\R-\{0\}\qquad  \lim_{z\in\mathbb{C}^+\to\lambda
} m_{F_{\widetilde
H, \wt c}}(z)\equiv
\breve{m}_{F_{\widetilde H, \wt c}}(\lambda)
\]
 exists, and that
$F_{\widetilde H, \wt c}$ has a continuous derivative
$F_{\widetilde H, \wt c}'=\pi^{-1}\im[\breve{m}_{F_{\widetilde H,
\wt c}} ]$ on
$(0,+\infty)$. In the case
$\wt c < 1$, $F_{\widetilde H, \wt c}$ has a continuous derivative on
all of $\R$
with $F^\prime_{\widetilde H, \wt c}
\equiv0$ on $(-\infty, 0]$.

For a grid $Q$ on the real line and for two c.d.f.'s $G_1$
and $G_2$, define
\[
\Vert G_1 - G_2\Vert_{Q} \equiv\sup_{t \in Q} |G_1(t) - G_2(t)|.
\]

The following theorem shows that both $F$ and $H$ can be estimated
consistently via an idealized algorithm.

\begin{theorem}\label{theo:cons1} Let $\{Q_n\}$ be a sequence of grids
on the real line eventually covering the support of $F$
with corresponding grid sizes $\{\gs_n\}$
satisfying $\gs_n \to0$.
Let $\{\widehat c_n\}$ be a sequence of positive real numbers with
$\wh c_n \to c$.
Let $\widehat H_n$ be defined as
%
\begin{equation}\label{e:min-global}
\widehat H_n \equiv\argmin_{\widetilde H} \Vert F_{\widetilde H, \wh c_n}
- F_n\Vert_{Q_n},
\end{equation}
where $\widetilde H$ is a probability measure.

Then we have
\emph{(i)} $F_{\widehat H_n, \wh c_n} \Rightarrow F $ a.s.; and
\emph{(ii)} $\widehat H_n \Rightarrow H$ a.s.
\end{theorem}

The algorithm used in the theorem is not practical for two
reasons. First, it is not possible to optimize over all
probability
measures $\widetilde H$. But similarly to\vadjust{\goodbreak} El Karoui (\citeyear{karoui:2008}), we can show
that it is sufficient to optimize over all probability measures that
are sums of atoms, the location of which is restricted to a fixed-size
grid, with the grid size vanishing asymptotically.

\begin{corollary}\label{cor:min-local1}
Let $\{Q_n\}$ be a sequence of grids
on the real line eventually covering the support of $F$
with corresponding grid sizes $\{\gamma_n\}$
satisfying $\gs_n \to0$.
Let $\{\widehat c_n\}$ be a sequence of positive real numbers with
$\wh c_n \to c$.
Let ${\mathcal P}_n$ denote the set of
all probability measures that are sums of atoms
belonging to the grid $\{J_n/T_n, (J_n + 1)/T_n, \ldots, K_n/T_n\}$
with $T_n \to\infty$,
$J_n$ being the largest integer satisfying $J_n/T_n \le\lambda_1$, and
$K_n$ being the smallest integer\vspace*{1pt} satisfying $K_n/T_n \ge\lambda_p$.
Let $\widehat H_n$ be defined as
%
\begin{equation}\label{e:min-local1}
\widehat H_n \equiv\argmin_{\widetilde H \in{\mathcal P}_n}
\Vert F_{\widetilde H, \wh c_n} - F_n\Vert_{Q_n}.
\end{equation}
Then we have
\emph{(i)} $F_{\widehat H_n, \wh c_n} \Rightarrow F $ a.s.; and
\emph{(ii)} $\widehat H_n \Rightarrow H$ a.s.
\end{corollary}

But even restricting the optimization over a manageable set of
probability measures is not quite practical yet for a second
reason. Namely, to compute $F_{\widetilde H, \wh c_n}$ exactly for a
given $\widetilde H$, one
would have to (numerically) solve the Mar\v{c}enko--Pastur equation
for an infinite
number of points. In practice, we can only afford to solve the
equation for a finite number of points and then approximate
$F_{\widetilde
H, \wh c_n}$ by trapezoidal integration. Fortunately, this
approximation does
not negatively affect the consistency of our estimators.\looseness=1

Let $G$ be a c.d.f. with continuous density $g$ and compact support
$[a, b]$.
For a grid $Q \equiv\{\ldots, t_{-1}, t_0, t_1, \ldots\}$
covering the support of $G$,
the approximation to $G$ via trapezoidal integration over the grid
$Q$, denoted by $\widehat G_Q$, is obtained as
follows. For $t \in[a, b]$, let $J_{lo} \equiv\max\{k\dvtx t_k \le a\}$
and $J_{hi} \equiv\min\{k\dvtx t < t_k\}$. Then
%
\begin{equation}\label{e:trapez1}
\widehat G_Q(t) \equiv
\sum_{k=J_{lo}}^{J_{hi}-1} \frac{(t_{k+1} - t_k) [g(t_k) +
g(t_{k+1})]}{2}.
\end{equation}

Now turn to the special case $G \equiv F_{\widetilde H, \wt c}$ and $Q
\equiv Q_n$. In
this case, we denote the approximation to $F_{\widetilde
H, \wt c}$ via trapezoidal integration over the grid $Q_n$ by $\widehat
F_{\widetilde H, \wt c; Q_n}$.

\begin{corollary}\label{cor:min-local2}
Assume the same assumptions as in Corollary~\ref{cor:min-local1}.
Let~$\widehat H_n$ be defined as
%
\begin{equation}\label{e:min-local2}
\widehat H_n \equiv\argmin_{\widetilde H \in{\mathcal P}_n}
\Vert\widehat F_{\widetilde H, \wh c_n; Q_n} - F_n\Vert_{Q_n}.
\end{equation}
Let $\breve m_{F_{\widehat H_n, \wh c_n}}(\lambda)$, $\widehat S_n$,
and $\widehat P_n$
be defined as in \eqref{e:m-breve-hat}, \eqref{e:nl-estimator}
and \eqref{e:nl-estimator-precision},
respectively.
Then:
\begin{longlist}[(iii)]
\item[(i)] $F_{\widehat H_n, \wh c_n} \Rightarrow F $ a.s.
\item[(ii)] $\widehat H_n \Rightarrow H$ a.s.
\item[(iii)] For any $\wt\delta> 0$, $\breve m_{F_{\widehat H_n, \wh
c_n}}(\lambda)\to
\breve{m}_F(\lambda)$ a.s. uniformly in $\lambda\in[\wt z_1 + \wt
\delta,
\wt z_2 - \wt\delta]$.
\item[(iv)] $\Vert\widehat S_n - S_n^{or}\Vert \to0$ a.s.
\item[(v)] $\Vert\widehat P_n - P_n^{or}\Vert \to0$ a.s.
\end{longlist}
\end{corollary}

\subsection{Implementation details}\label{sec:procedure}

\subsubsection*{Decomposition of the c.d.f. of population eigenvalues}
As discussed before,
it is not practical to search over the set of all possible c.d.f.'s
$\widetilde{H}$. Following \citet{karoui:2008}, we project $H$ onto a
certain basis of c.d.f.'s $(M_k)_{k=1,\ldots,K}$, where $K$ goes to
infinity along
with $n$ and $p$. The projection of $H$ onto this basis is
given by the nonnegative weights $w_1,\ldots,w_K$, where
%
\begin{equation}
\label{eq:weights}
\forall t\in\R\qquad H(t) \approx\widetilde{H} (t) \equiv\sum
_{k=1}^Kw_kM_k(t) \quad\mbox{and}\quad \sum_{k=1}^Kw_k=1.
\end{equation}
Thus, our estimator for $F$ will be a solution to the
Mar\v{c}enko--Pastur equation for~$\widetilde{H}$ given by equation
(\ref{eq:weights}) for some $(w_k)_{k=1,\ldots,K}$, and for
$\widetilde c \equiv p/n$.
It is just a matter of searching over all sets of nonnegative weights
summing up to one.

\subsubsection*{Choice of basis}
We base the c.d.f.'s $(M_k)_{k=1,\ldots,K}$ on a grid of $p$ equally
spaced points on the interval $[\lambda_1,\lambda_p]$.
%
\begin{equation}\label{e:basis}
x_i \equiv\lambda_1+\frac{i-1}{p}(\lambda_p-\lambda_1) \qquad\mbox{for
} i=1,\ldots,p.
\end{equation}
Thus $x_1=\lambda_1$ and $x_p=\lambda_p$. We then form the basis $\{
M_1,\ldots,M_k\}$ as the union of three families of c.d.f.'s:
\begin{enumerate}[(3)]
\item[(1)] the indicator functions $\bone_{[x_i,+\infty)}$
($i=1,\ldots,p$);
\item[(2)] the c.d.f.'s whose derivatives are linearly increasing on
the interval $[x_{i-1},x_i]$ and zero everywhere else ($i=2,\ldots,p$);
\item[(3)] the c.d.f.'s whose derivatives are linearly decreasing on
the interval $[x_{i-1},x_i]$ and zero everywhere else ($i=2,\ldots,p$).
\end{enumerate}
This list yields a basis $(M_k)_{k=1,\ldots,K}$ of dimension
$K=3p-2$. Notice that by the theoretical results of
Section~\ref{ss:H-cons}, it would be sufficient to use the first
family only. Including the second and third families in addition cannot
make the
approximation to $H$ any worse.

\subsubsection*{Trapezoidal integration}
For a given $\widetilde{H} \equiv\sum_{k=1}^Kw_kM_k$, it is computationally
too expensive (in the context of an optimization procedure) to
solve
the Mar\v{c}enko--Pastur\vadjust{\goodbreak} equation for $m_F(z)$ over all $z\in\C^+$. It
is more efficient to solve the Mar\v{c}enko--Pastur equation only for
$\breve{m}_F(x_i)$ $(i=1,\ldots,p)$, and to use the trapezoidal
approximation formula to deduce from it $F(x_i)$
$(i=1,\ldots,p)$. The trapezoidal rule gives
%
\begin{eqnarray}\label{eq:trapeze}
\forall i=1,\ldots,p\qquad F(x_i)&=&\sum_{j=1}^{i-1}\frac
{x_{j+1}-x_{j-1}}{2} F'(x_j)
+\frac{x_{i}-x_{i-1}}{2} F'(x_i)\nonumber
\\
&=&\sum_{j=1}^{i-1}\frac{ (x_{j+1}-x_{j-1} )\im[\breve{m}_F(x_j)
]}{2\pi}
\\
&&{}+\frac{ (x_{i}-x_{i-1} )\im[\breve{m}_F(x_i) ]}{2\pi},\nonumber
\end{eqnarray}
with the convention $x_{0} \equiv0$.

\subsubsection*{Objective function}
The objective function measures the distance between $F_n$ and the $F$
that solves the Mar\v{c}enko--Pastur equation for $\widetilde{H}
\equiv
\sum_{k=1}^Kw_kM_k$ and for $\widetilde c \equiv p/n$.
Traditionally, $F_n$ is defined as
c\`{a}dl\`{a}g, that is, $F_n(\lambda_1)=1/p$ and
$F_n(\lambda_p)=1$. However, there is a certain degree of
arbitrariness in this convention: why is $F_n(\lambda_p)$ equal to one
but $F_n(\lambda_1)$ not equal to zero? By symmetry, there is no
a priori justification for specifying that the largest eigenvalue
is closer to the supremum of the support of $F$ than the smallest to
its infimum. Therefore, a~different convention might be more
appropriate in this case, which leads us to the following definition:
%
\begin{equation}
\label{eq:convention}
\forall i=1,\ldots,p\qquad \widehat{F}_n(\lambda_i) \equiv\frac
{i}{p}-\frac{1}{2p}.
\end{equation}
This choice restores a certain element of symmetry to the treatment of
the smallest vs. the largest eigenvalue.
From equation (\ref{eq:convention}), we deduce $\widehat{F}_n(x_i)$, for
$i=2,\ldots,p-1$, by linear interpolation.
With a sup-norm error penalty, this leads to the following objective function:
%
\begin{equation}
\label{eq:objective}
\max_{i=1,\ldots,p} |F(x_i)-\widehat{F}_n(x_i) |,
\end{equation}
where $F(x_i)$ is given by equation (\ref{eq:trapeze}) for $i=1,\ldots
,p$. Using equation (\ref{eq:trapeze}), we can rewrite this objective
function as
\[
\max_{i=1,\ldots,p} \Biggl|\sum_{j=1}^{i-1}\frac{ (x_{j+1}-x_{j-1} )\im
[\breve{m}_F(x_j) ]}{2\pi}
+\frac{ (x_{i}-x_{i-1} )\im[\breve{m}_F(x_i) ]}{2\pi}-\widehat
{F}_n(x_i) \Biggr|.
\]

\subsubsection*{Optimization program}
We now have all the ingredients needed to state the optimization
program that will extract the estimator of $\breve{m}_F(x_1),\ldots
,\allowbreak\breve{m}_F(x_p)$ from the observations $\lambda_1,\ldots, \lambda
_p$. It is the following:
\[
\mathop{\min_{m_1,\ldots,m_p}}_{w_1,\ldots,w_K}\max
_{i=1,\ldots,p} \Biggl|\sum_{j=1}^{i-1}\frac{ (x_{j+1}-x_{j-1} )\im[m_j
]}{2\pi}
+\frac{ (x_{i}-x_{i-1} )\im[m_i ]}{2\pi}-\widehat{F}_n(x_i) \Biggr|\vadjust{\goodbreak}
\]
subject to
%
\begin{eqnarray}\label{eq:mpconstraint}
\forall j=1,\ldots,p\qquad m_j&=&\sum_{k=1}^K\int_{-\infty}^{+\infty
}\frac{w_k}{t [1-(p/n)-(p/n) x_j m_j ]-x_j}\, dM_k(t),\nonumber\hspace*{-35pt}
\\
\sum_{k=1}^Kw_k&=&1,\nonumber
\\[-12pt]
\\[-4pt]
\forall j=1,\ldots,p\qquad m_j&\in&\C^+,\nonumber\hspace*{-35pt}
\\
\forall k=1,\ldots,K\qquad w_k&\geq&0. \nonumber\hspace*{-35pt}
\end{eqnarray}
The key is to introduce the variables
$m_j\equiv\breve{m}_F(x_j)$, for $j=1,\ldots,p$. The constraint
in equation (\ref{eq:mpconstraint}) imposes that $m_j$ is the solution
to the Mar\v{c}enko--Pastur equation evaluated as $z\in\C^+\to x_j$ when
$ \widetilde{H} =\sum_{k=1}^Kw_kM_k$.

\subsubsection*{Real optimization program}
In practice, most optimizers only accept real variables. Therefore it
is necessary to decompose $m_j$ into its real and imaginary parts:
$a_j \equiv\re[m_j]$ and $b_j \equiv\im[m_j]$. Then we can optimize
separately over the two sets of real variables $a_j$ and $b_j$ for
$j=1,\ldots,p$. The Mar\v{c}enko--Pastur constraint in equation (\ref
{eq:mpconstraint})
splits into two constraints: one for the real part and the other for
the imaginary part. The reformulated optimization program is
%
\begin{equation}
\label{eq:program}\qquad
\mathop{\mathop{\min_{a_1,\ldots,a_p}}_{b_1,\ldots,b_p}}_{
w_1,\ldots,w_K}
\max_{i=1,\ldots,p} \Biggl|\sum_{j=1}^{i-1}\frac{ (x_{j+1}-x_{j-1}
)b_j}{2\pi}
+\frac{ (x_{i}-x_{i-1} )b_i}{2\pi}-\widehat{F}_n(x_i) \Biggr|
\end{equation}
subject to
%
\begin{eqnarray}
\label{eq:MPreal}%
&&\forall j=1,\ldots,p\nonumber\hspace*{-35pt}
\\[-8pt]
\\[-8pt]
&&\quad a_j= \sum_{k=1}^K\int_{-\infty}^{+\infty}\re
\biggl\{\frac{w_k}{t [1-(p/n)-(p/n) x_j(a_j+ib_j) ]-x_j} \biggr\}\,dM_k(t),\nonumber\hspace*{-35pt}
\\
\label{eq:MPimag}&&\forall j=1,\ldots,p
\nonumber\hspace*{-35pt}
\\[-8pt]
\\[-8pt]
&&\quad b_j= \sum_{k=1}^K\int_{-\infty}^{+\infty}\im
\biggl\{\frac{w_k}{t [1-(p/n)-(p/n) x_j(a_j+ib_j) ]-x_j} \biggr\}\,dM_k(t),\nonumber\hspace*{-35pt}
\\
&&\sum_{k=1}^Kw_k=1,\hspace*{-35pt}
\\
&&\forall j=1,\ldots,p\qquad b_j\geq0,\hspace*{-35pt}
\\
\label{eq:nonnegative}&&\forall k=1,\ldots,K\qquad w_k\geq0.\hspace*{-35pt}
\end{eqnarray}

\begin{remark}
Since the theory
of Sections~\ref{sec:m} and~\ref{ss:H-cons} partly
assumes that $m_j$ belongs to a compact set in $\C^+$ bounded away
from the real line, we might want to add to the real optimization program
the constraints that $-1/\varepsilon\leq a_j\leq
1/\varepsilon$ and that $\varepsilon\leq b_j \leq1/\varepsilon$,
for some small
$\varepsilon>0$. Our simulations indicate that for a small value of
$\varepsilon$ such as $\varepsilon=10^{-6}$, this makes no difference
in practice.
\end{remark}

\subsubsection*{Sequential linear programming}
While the optimization program defined in equations
(\ref{eq:program})--(\ref{eq:nonnegative}) may appear daunting at first
sight because of its non-convexity,
it is, in fact, solved quickly and efficiently by off-the-shelf
optimization software implementing Sequential Linear Programming
(SLP). The key is to linearize equations
(\ref{eq:MPreal})--(\ref{eq:MPimag}), the two constraints that embody
the Mar\v{c}enko--Pastur equation, around an approximate solution
point. Once they are
linearized, the optimization program
(\ref{eq:program})--(\ref{eq:nonnegative}) becomes a standard Linear
Programming (LP) problem, which can be solved very quickly. Then we
linearize again equations (\ref{eq:MPreal})--(\ref{eq:MPimag}) around
the new point, and this generates a new LP problem; hence the name:
\emph{Sequential} Linear Programming. The software iterates until a
satisfactory degree of convergence is achieved. All of this is handled
automatically by the SLP optimizer. The user only needs to specify the
problem (\ref{eq:program})--(\ref{eq:nonnegative}), as well as
some starting point, and then launch the SLP optimizer. For our
SLP optimizer, we selected a standard off-the-shelf commercial
software: SNOPT\texttrademark\
Version~7.2--5; see \citet{gill:murray:saunders:2002}.
While SNOPT\texttrademark\ was originally designed for sequential
quadratic programming, it also handles SLP, since linear programming
can be viewed as a particular case of quadratic programming
with no quadratic term.

\subsubsection*{Starting point}
A neutral way to choose the starting point is to place equal weights
on all the c.d.f.'s in our basis: $w_k \equiv1/K
(k=1,\ldots,K)$. Then it is necessary to solve the Mar\v
{c}enko--Pastur equation
numerically once \emph{before} launching the SLP optimizer, in order to
compute the values of $\breve{m}_F(x_j)$ $(j=1,\ldots,p)$ that
correspond to this initial choice of $\widetilde{H}=\sum
_{k=1}^KM_k/K$. The
initial values for $a_j$ are taken to be
$\re[\breve{m}_F(x_j) ]$, and
$\im[\breve{m}_F(x_j) ]$ for $b_j$
$(j=1,\ldots,p)$. If the choice of equal weights $w_k \equiv1/K$ for
the starting point does not lead to convergence of the optimization
program within a~pre-specified limit on the maximum number of iterations,
we choose random weights $w_k$ generated i.i.d. $\sim\operatorname{Uniform}[0,1]$ (rescaled to sum up to one), repeating
this process until
convergence finally occurs. In the vast majority of
cases, the optimization program already converges on the first
try. For example, over 1000 Monte Carlo simulations using the design
of Section~\ref{sub:MC_convergence} with $p = 100$ and $n = 300$,
the optimization program converged on the first try 994 times and on
the second try the remaining 6 times.
%
\begin{figure}

\includegraphics{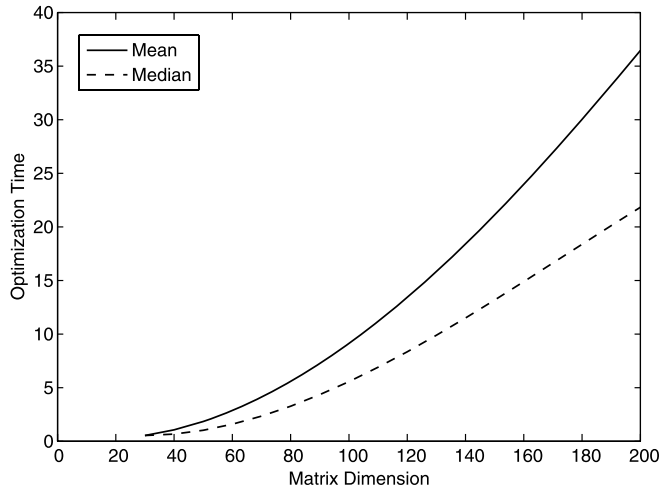}

\caption{Mean and median CPU times (in seconds) for optimization program as function of
matrix dimension. The design is the one of
Section~\protect\ref{sub:MC_convergence} with $n = 3 p$.
Every point is the result of $1000$ Monte Carlo simulations.}\label{fig:timer}
\end{figure}

\subsubsection*{Optimization time}

Figure~\ref{fig:timer} gives some information on how the
optimization time increases with the matrix
dimension.

The main reason for the rate at which the optimization time increases
with $p$ is that
the number of grid points in \eqref{e:basis} increases linearly in
$p$. This linear rate is not a requirement for our asymptotic
results. Therefore, if necessary, it is possible to pick a
less-than-linear rate of increase in the number of grid points to
speed up the optimization for very large matrices.

\subsubsection*{Estimating the covariance matrix} Once the SLP optimizer
has converged, it generates optimal values $(a_1^*,\ldots,a_p^*)$,
$(b_1^*,\ldots,b_p^*)$ and $(w_1^*,\ldots,w_K^*)$. The first two sets
of variables at the optimum are used to estimate the oracle
shrinkage factors.
From the
reconstructed $\breve{m}_F^*(x_j)\equiv a_j^*+i b_j^*$, we deduce by linear
interpolation $\breve{m}_F^*(\lambda_j)$, for $j=1,\ldots,p$. Our
estimator of the covariance matrix $\widehat{S}_n$ is built
by keeping the same eigenvectors as the sample covariance matrix, and
dividing each sample eigenvalue $\lambda_j$ by the following
correction factor:

%
\[
\biggl|1-\frac{p}{n}-\frac{p}{n} \lambda_j \breve{m}_F^*(\lambda_j) \biggr|^2.
\]
%
Corollary~\ref{cor:min-local2} assures us that the resulting \textit{bona
fide} nonlinear shrinkage estimator is asymptotically equivalent to
the oracle estimator $S_n^{or}$.
Also, we can see that, as the concentration $\widehat{c}_n= p/n$ gets closer
to zero, that is, as we get closer to fixed-dimension asymptotics, the
magnitude of the correction becomes smaller. This makes
sense because under fixed-dimension asymptotics the sample covariance
matrix is a consistent estimator of the population covariance matrix.

\subsubsection*{Estimating the precision matrix} The output of the
same optimization process can also be used to estimate the oracle
shrinkage factors for the precision matrix.
Our estimator of the precision
matrix $\Sigma_n^{-1}$ is built by keeping the same eigenvectors as
the sample
covariance matrix, and multiplying the inverse $\lambda_j^{-1}$ of
each sample eigenvalue by the following correction factor:
\[
1-\frac{p}{n}-2 \frac{p}{n} \lambda_j
\re[\breve{m}_F^*(\lambda_j) ].
\]
Corollary~\ref{cor:min-local2} assures us that the resulting \textit{bona
fide} nonlinear shrinkage estimator is asymptotically equivalent to
the oracle estimator $P_n^{or}$.

\subsubsection*{Estimating $H$} We point out that the optimal values
$(w_1^*,\ldots,w_K^*)$ generated from the SLP optimizer yield a
consistent estimate of $H$ in the following fashion:
\[
H^* \equiv\sum_{k=1}^K w_k^* M_k.
\]

This estimator could be considered an alternative to the estimator
introduced by \citet{karoui:2008}.
The most salient difference between the two optimization
algorithms is that our objective function tries to match $F_n$ on~$\R$,
whereas his objective function tries to match (a function of) $m_{F_n}$
on~$\C^+$. The deeper we go into~$\C^+$, the more ``smoothed-out'' is the
Stieltjes transform, as it is an analytic function; therefore, the
more information is lost. However, the approach of \citet{karoui:2008}
cannot get too close to the real line because $m_{F_n}$ starts looking
like a sum of Dirac functions (which are very ill-behaved) as one gets
close to the real line, since $F_n$ is a step
function. In a sense, the approach of \citet{karoui:2008} is to match
a smoothed-out version of a sum of ill-behaved Diracs. In this
situation, knowing ``how much to smooth'' is rather delicate, and even
if it is done well, it still loses information. By contrast, we have
no information loss because we operate directly on the real line, and
we have no problems with Diracs because we match $F_n$ instead of its
derivative. The price to pay is that our optimization program is not
convex, whereas the one of \citet{karoui:2008} is. But extensive
simulations reported in
the next section show that off-the-shelf nonconvex optimization
software---as the commercial package SNOPT---can handle this
particular type of a nonconvex problem in
a fast, robust and efficient manner.

It would have been of additional interest to compare our estimator
of $H$ to the one of \citet{karoui:2008} in some
simulations. But when we tried to implement his estimator according to
the implementation details provided,
we were not able to match the results presented in his
paper. Furthermore, we were not able\vadjust{\goodbreak} to obtain his original software.
As a result, we
cannot make any definite statements concerning the performance of
our estimator of $H$ compared to the one of \citet{karoui:2008}.

\begin{remark}[(Cross-validation estimator)]\label{r:cv}
The implementation of our nonlinear shrinkage estimators is
not trivial and also requires the use of a~third-party SLP
optimizer. It is therefore of interest whether an alternative version
exists that is easier to implement and exhibits (nearly) as good
finite-sample properties.

To this end an anonymous referee suggested to estimate the quantities~$d_i^*$ of \eqref{e:star} by a leave-one-out cross-validation
method. In particular, let $(\lambda_i[k], \ldots,\allowbreak \lambda_p[k]);
(u_1[k], \ldots, u_p[k])$ denote a system of eigenvalues and
eigenvectors of the sample covariance matrix computed from all the
observed data, except for the $k$th observation. Then $d_i^*$
of \eqref{e:star} can be approximated by
\[
d_i^{cv} \equiv\frac{1}{n} \sum_{k=1}^n (u_i[k]^\prime y_k)^2,
\]
where the $p \times1$ vector $y_k$ denotes the $k$th row of the matrix
$Y_n \equiv X_n
\Sigma_n^{1/2}$.

The motivation here is that
\[
(u_i[k]^\prime y_k)^2 = u_i[k]^\prime y_k y_k^\prime u_i[k],
\]
where $y_k$ is independent of $u_i[k]$ and $\E(y_k y_k^\prime) =
\Sigma_n$ (even though $y_k y_k^\prime$ is of
rank one only).

We are grateful for this suggestion, since the cross-validation
quantities~$d_i^{cv}$ can be computed without the use of any
third-party optimization software, and the corresponding computer code
is very short.

On the other hand, the cross-validation estimator has three disadvantages.
First, when $p$ is large, it takes
much longer to compute the cross-validation estimator. The reason is
that the spectral decomposition of a $p \times p$ covariance matrix has to
be computed $n$ times as opposed to only one time. Second,
the cross-validation method only applies to the
estimation of the covariance matrix $\Sigma_n$ itself. It is not clear
how to adapt this method to the (direct) estimation of the
precision matrix $\Sigma_n^{-1}$ or any other smooth function of
$\Sigma_n$. Third, the performance of the cross-validation estimator
cannot match the performance of our method; see Section~\ref{ss:additional}.
\end{remark}

\begin{remark}
\label{rm:mestre}
Another approach proposed recently is the one of
\citet{mestre:lagunas:2006}. They use so-called ``G-estimation,''
that is, asymptotic results that assume the sample size $n$ and the
matrix dimension
$p$ go to infinity together, to derive minimum variance
beam formers in the context of the spatial filtering of electronic
signals. There are several differences between their paper\vadjust{\goodbreak} and the
present one. First, \citet{mestre:lagunas:2006} are interested in an
optimal $p \times1$ weight vector $w_{opt}$ given by
\[
w_{opt} \equiv\argmin_{w} w^\prime\Sigma_n w\qquad  \mbox{subject to } w^\prime s_d = 1,
\]
where $s_d$ is a $p \times1 $ vector containing signal
information. Consequently, \citet{mestre:lagunas:2006} are ``only''
interested in a certain functional of $\Sigma_n$, while
we are interested in the full covariance matrix $\Sigma_n$
and also in the full precision matrix $\Sigma_n^{-1}$. Second, they
use the real
Stieltjes transform, which is different from the more conventional
complex Stieltjes transform used in random matrix theory and in the
present paper. Third, their random variables are complex whereas ours
are real. The cumulative impact of these differences is best
exemplified by the
estimation of the precision matrix: \citeauthor{mestre:lagunas:2006} [(\citeyear{mestre:lagunas:2006}),
page 76] recommend $(1-p/n) S_n^{-1}$, which is just a rescaling of the
inverse of the sample covariance matrix, whereas our
Section~\ref{ss:oracle:precision} points to a highly nonlinear
transformation
of the eigenvalues of the sample covariance matrix.
\end{remark}

\section{Monte Carlo simulations}
\label{sec:MC}

In this section, we present the results of various sets of Monte Carlo\vspace*{1pt}
simulations designed to illustrate the finite-sample properties of the
nonlinear shrinkage estimator $\widehat S_n$. As detailed in
Section~\ref{sec:oracle}, the finite-sample optimal estimator in the
class of rotation-equivariant estimators is given by $S_n^*$ as
defined in \eqref{e:optimal-fs}.
Thus, the improvement of the shrinkage estimator $\widehat{S}_n$ over
the sample covariance matrix will be measured by how closely this
estimator approximates $S^*_n$ relative to the sample covariance
matrix. More specifically, we report the Percentage Relative
Improvement in Average Loss (PRIAL), which is defined as
%
\begin{equation}\label{e:prial}
\mathrm{PRIAL} \equiv\mathrm{PRIAL}(\widehat\Sigma_n) \equiv
100\times
\biggl\{1-\frac{\E[ \Vert\widehat{\Sigma}_n-S^*_n \Vert^2 ]}
{\E[ \Vert S_n-S^*_n \Vert^2 ]} \biggr\} \%,
\end{equation}
where $\widehat\Sigma_n$ is an arbitrary estimator of $\Sigma_n$.
By definition, the PRIAL of $S_n$ is 0\%, while the PRIAL of $S_n^*$ is
100\%.

Most of the simulations will be designed around a population
covariance matrix $\Sigma_n$ that has $20\%$ of its eigenvalues equal
to $1$, $40\%$ equal to $3$ and $40\%$ equal to~$10$. This is a
particularly interesting and difficult example introduced and analyzed
in detail by \citet{bai:silverstein:1998}. For concentration
values such as $c=1/3$ and below, it displays ``spectral separation;''
that is, the support of the distribution of sample eigenvalues is the
union of three disjoint intervals, each one corresponding to a Dirac of
population eigenvalues. Detecting this pattern and handling it
correctly is a real challenge for any covariance matrix estimation
method.

\subsection{Convergence}\label{sub:MC_convergence}

\begin{figure}

\includegraphics{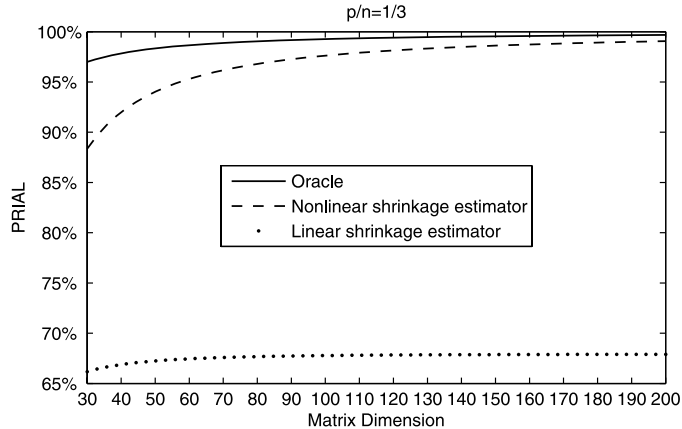}

\caption{Comparison of the nonlinear vs. linear shrinkage estimators. $20\%$ of
population eigenvalues are equal to $1$, $40\%$ are equal to $3$ and
$40\%$ are equal to $10$. Every point is the result of $1000$
Monte Carlo simulations.}\label{fig:convergence}
\end{figure}

The first set of Monte Carlo simulations shows how the nonlinear
shrinkage estimator $\widehat{S}_n$ behaves as the matrix dimension
$p$ and the sample size $n$ go to infinity together. We assume that
the concentration ratio\vadjust{\goodbreak} $\widehat{c}_n=p/n$ remains constant and equal
to $1/3$. For every value of $p$ (and hence $n$), we run 1000
simulations with normally distributed variables. The PRIAL is plotted
in Figure~\ref{fig:convergence}. For the sake of comparison, we also
report the PRIALs of the oracle $S_n^{or}$ and the optimal linear
shrinkage estimator
$\overline{S}_n$ developed by \citet{ledoit:wolf:2004a}.\vspace*{1pt}

%

\begin{figure}

\includegraphics{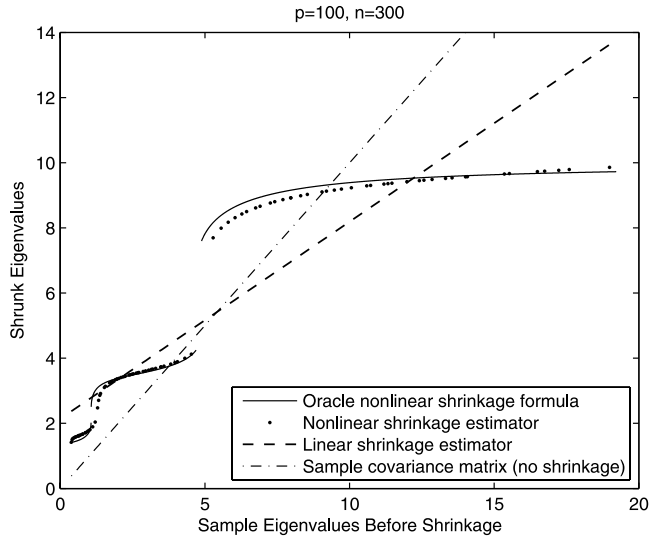}

\caption{Nonlinearity of the oracle shrinkage formula. $20\%$ of population
eigenvalues are equal to~$1$, $40\%$ are equal to $3$ and $40\%$ are
equal to $10$. $p=100$ and $n=300$.}\label{fig:nonlinearity}
\end{figure}

One can see that the performance of the nonlinear shrinkage
estimator~$\widehat S_n$ converges
quickly toward that of the oracle and of $S^*_n$.
Even for relatively small matrices of
dimension $p=30$, it realizes $88\%$ of the possible gains over the
sample covariance matrix. The optimal linear shrinkage estimator~$\overline{S}_n$ also performs well relative to the sample covariance
matrix, but the improvement is limited: in general, it does not
converge to $100\%$ under large-dimensional asymptotics. This is
because there are strong nonlinear effects in the optimal shrinkage of
sample eigenvalues. These effects are clearly visible in
Figure~\ref{fig:nonlinearity}, which plots a typical simulation result for
$p=100$.

%

One can see that the nonlinear shrinkage estimator $\widehat{S}_n$
shrinks the eigenvalues of the sample covariance matrix almost as if
it ``knew'' the correct shape of the distribution of population
eigenvalues. In particular,
the various curves and gaps of the oracle nonlinear
shrinkage formula are well picked up and followed by this
estimator. By contrast, the linear shrinkage estimator can only use
the best linear approximation to this highly nonlinear
transformation. We also plot the $45$-degrees line as a visual
reference to show what would happen if no shrinkage was applied to the
sample eigenvalues, that is, if we simply used $S_n$.

\subsection{Concentration}
\label{sub:MC_concentration}

\begin{figure}[b]

\includegraphics{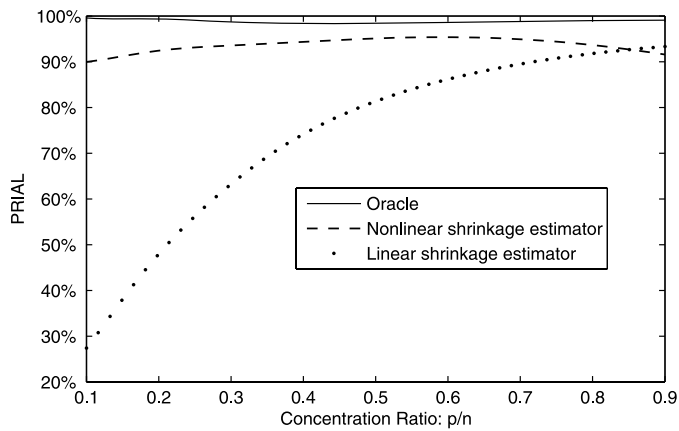}

\caption{Effect of varying the concentration ratio $\widehat{c}_n=p/n$. $20\%$
of population eigenvalues are equal to $1$, $40\%$ are equal to $3$
and $40\%$ are equal to $10$. Every point is the result of $1000$
Monte Carlo simulations.}\label{fig:concentration}
\end{figure}

The next set of Monte Carlo simulations shows how the PRIAL of the
shrinkage estimators varies as a function of the concentration ratio
$\widehat{c}_n=p/n$ if we keep the product $p\times n$ constant and
equal to $9000$. We keep the same population covariance matrix
$\Sigma_n$ as in Section~\ref{sub:MC_convergence}. For every value of
$p/n$, we run $1000$ simulations with normally distributed
variables. The respective PRIALs of $S_n^{or}$, $\widehat{S}_n$ and
$\overline{S}_n$ are plotted in Figure~\ref{fig:concentration}.

%

One can see that the nonlinear shrinkage estimator performs well across
the board, closely in line with the oracle, and always achieves at
least $90\%$ of the possible
improvement over the sample covariance matrix. By contrast, the linear
shrinkage estimator achieves relatively little improvement over the
sample covariance matrix when the concentration is low. This is
because, when the sample size is large relative to the matrix
dimension, there is a lot of precise information about the optimal
nonlinear way to shrink the sample eigenvalues that is waiting to be
extracted by a suitable nonlinear procedure. By contrast, when the
sample size is not so large, the information about the population
covariance matrix is relatively fuzzy; therefore a simple linear
approximation can achieve up to $93\%$ of the potential gains.\vspace*{-3pt}

\subsection{Dispersion}
\label{sub:MC_dispersion}

\begin{figure}[b]
\vspace*{-3pt}
\includegraphics{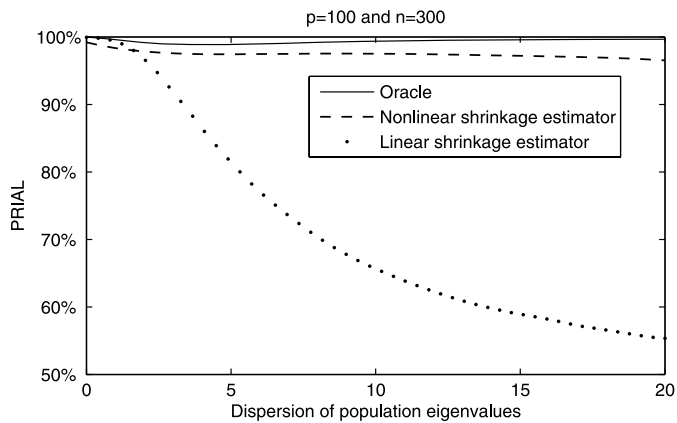}

\caption{Effect of varying the dispersion of population eigenvalues. $20\%$ of
population eigenvalues are equal to $1$, $40\%$ equal to $1+2 d/9$
and $40\%$ equal to $1+d$, where the dispersion parameter $d$ varies
from $0$ to $20$. $p=100$ and $n=300$. Every point is the result of
$1000$ Monte Carlo simulations.}\label{fig:dispersion}
\end{figure}

The third set of Monte Carlo simulations shows how the PRIAL of the
shrinkage estimators varies as a function of the dispersion of
population eigenvalues. We take a population covariance matrix
$\Sigma_n$ with $20\%$ of its eigenvalues equal to $1$, $40\%$ equal
to $1+2 d/9$ and $40\%$ equal to $1+d$, where the dispersion
parameter $d$ varies from $0$ to $20$. Thus, for $d=0$, $\Sigma_n$ is
the identity matrix and, for $d=9$, $\Sigma_n$ is the same matrix as
in Section~\ref{sub:MC_convergence}. The sample size is $n=300$ and
the matrix dimension is $p=100$. For every value of $d$, we run $1000$
simulations with normally distributed variables. The respective PRIALs
of $S_n^{or}$, $\widehat{S}_n$ and $\overline{S}_n$ are plotted in
Figure~\ref{fig:dispersion}.

%

One can see that the linear shrinkage estimator $\overline{S}_n$ beats
the nonlinear shrinkage estimator $\widehat{S}_n$ for very low
dispersion levels. For example, when $d=0$, that is,
when the population covariance matrix is equal to the identity matrix,
$\overline{S}_n$ realizes $99.9\%$ of the possible improvement over
the sample covariance matrix, while $\widehat{S}_n$ realizes ``only''
$99.4\%$ of the possible improvement. This is because, in this case,
linear shrinkage is optimal or (when $d$ is strictly positive but
still small) nearly optimal Hence there is nothing too little to be
gained by\vadjust{\goodbreak}
resorting to a nonlinear shrinkage method.
However, as dispersion increases, linear shrinkage delivers less and
less improvement over the sample covariance matrix, while nonlinear
shrinkage retains a PRIAL above $96\%$, and close to that of the oracle.

\subsection{Fat tails}

We also have some results on the effect of non-normality on the
performance of the shrinkage estimators. We take the same population
covariance matrix as in Section~\ref{sub:MC_convergence}, that is,
$\Sigma_n$ has $20\%$ of its eigenvalues equal to $1$, $40\%$ equal to
$3$ and $40\%$ equal to $10$. The sample size is $n=300$, and the matrix
dimension is $p=100$. We compare two types of random variates: a
Student~$t$ distribution with $\mathrm{df}=3$ degrees of freedom, and a
Student~$t$ distribution with $\mathrm{df}=\infty$ degrees of freedom
(which is the Gaussian distribution). For each number of degrees of
freedom $\mathrm{df}$, we run $1000$ simulations. The respective PRIALs
of $S_n^{or}$, $\widehat{S}_n$ and $\overline{S}_n$ are summarized in
Table~\ref{table:fat_tails}.

%
\begin{table}
\tablewidth=\textwidth
\tabcolsep=0pt
\caption{Effect of nonnormality. $20\%$ of population eigenvalues are equal to
$1$, $40\%$ are equal to $3$ and $40\%$ are equal to $10$. $1000$
Monte Carlo simulations with $p=100$ and $n=300$}\label{table:fat_tails}
  \begin{tabular*}{\textwidth}{@{\extracolsep{\fill}}lcccc@{}}
\hline
& \multicolumn{2}{c}{\textbf{Average squared}} & \multicolumn{2}{c}{\multirow{2}{*}{\textbf{PRIAL}}}
\\
& \multicolumn{2}{c}{\textbf{Frobenius loss}} &
\\[-5pt]
& \multicolumn{2}{c}{\hrulefill} & \multicolumn
{2}{c@{}}{\hrulefill}
\\
& $ \mathbf{df}\boldsymbol{=3} $ & $\mathbf{df}\boldsymbol{=\infty}$ & $ \mathbf{df}\boldsymbol{=3} $ &
$\mathbf{df}\boldsymbol{=\infty}$
\\
\hline
Sample covariance matrix & 5.856 & 5.837 & 0\%\phantom{.0} & 0\%\phantom{.0}
\\
Linear shrinkage estimator & 1.883 & 1.883 & 67.84\% & 67.74\%
\\
Nonlinear shrinkage estimator & 0.128 & 0.133 & 97.81\% & 97.71\%
\\
Oracle & 0.043 & 0.041 & 99.27\% & 99.30\%
\\ \hline
\end{tabular*}
\end{table}

One can see that departure from normality does not have
any noticeable effect on performance.

\subsection{Precision matrix}

The next set of Monte Carlo simulations focuses on estimating the
precision matrix $\Sigma_n^{-1}$. The definition of the PRIAL, in
this subsection only, is given by
%
\begin{equation}\label{e:prial-inv}
\mathrm{PRIAL} \equiv\mathrm{PRIAL}(\widehat\Pi_n) \equiv
100\times
\biggl\{1-\frac{\E[ \Vert\widehat{\Pi}_n-P^*_n \Vert^2 ]}
{\E[ \Vert S_n^{-1}-P^*_n \Vert^2 ]} \biggr\} \%,
\end{equation}
where $\widehat\Pi_n$ is an arbitrary estimator of $\Sigma_n^{-1}$.
By definition, the PRIAL of $S_n^{-1}$ is 0\% while the PRIAL of $P_n^*$
is 100\%.

\begin{figure}

\includegraphics{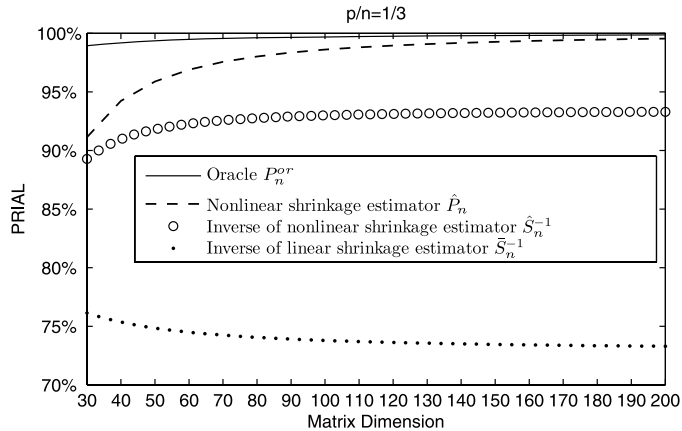}

\caption{Estimating the precision matrix. $20\%$ of population eigenvalues are
equal to $1$, $40\%$ are equal to $3$ and $40\%$ are equal to $10$.
Every point is the result of $1000$ Monte Carlo simulations.}\label{fig:inverse}
\end{figure}

We take the same
population eigenvalues as in Section~\ref{sub:MC_convergence}. The
concentration ratio $\widehat{c}_n=p/n$ is set to the value $1/3$. For
various values of $p$ between $30$ and $200$, we run 1000 simulations
with normally distributed variables.
The respective PRIALs of $P_n^{or}$, $\widehat{P}_n$, $\widehat{S}_n^{-1}$
and $\overline{S}_n^{-1}$
are plotted in Figure~\ref{fig:inverse}.

%

One can see that the nonlinear shrinkage method seems to be just as
effective for\vadjust{\goodbreak} the purpose of estimating the precision matrix as it is
for the purpose of estimating the covariance matrix itself. Moreover,
there is a clear benefit in directly estimating the precision matrix
by means of $\widehat P_n$ as opposed to the indirect estimation by
means of $\widehat S_n^{-1}$ (which on its own significantly
outperforms $\overline{S}_n^{-1}$).

\subsection{Shape}

\begin{figure}

\includegraphics{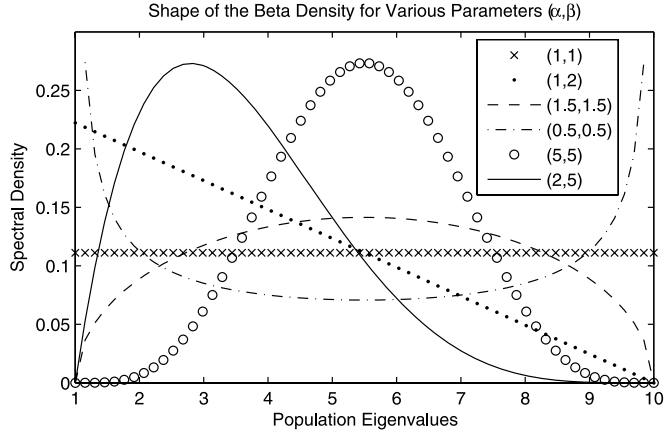}

\caption{Shape of the beta density for various parameter values. The support of
the beta density has been shifted to the interval $[1,10]$ by a linear
transformation. To enhance clarity, the densities corresponding to the
parameters $(2,1)$ and $(5,2)$ have been omitted, since they are
symmetric to $(1,2)$ and $(2,5)$, respectively, about the mid-point of
the support.}\label{fig:beta}
\end{figure}

Next, we study how the nonlinear shrinkage estimator $\widehat S_n$
performs for a wide variety of shapes of population spectral densities.
This requires using a family of distributions with bounded support and
which, for various parameter values, can take on different shapes. The
best-suited family for this purpose is the beta distribution.
The c.d.f. of the beta distribution with parameters $(\alpha,\beta)$ is
\[
\forall x\in[0,1]\qquad F_{(\alpha,\beta)}(x)=\frac{\Gamma(\alpha+\beta
)}{\Gamma(\alpha)\Gamma(\beta)}\int_0^xt^{\alpha-1}(1-t)^{\beta-1}\,dt.
\]
While the support of the beta distribution is $[0,1]$, we shift it to
the interval $[1,10]$ by applying a linear transformation. Thanks to
the flexibility of the beta family of densities, selecting different
parameters $(\alpha,\beta)$ enables us to generate eight different shapes
for the population spectral density: rectangular $(1,1)$, linearly
decreasing triangle $(1,2)$, linearly increasing triangle $(2,1)$,
circular $(1.5,1.5)$, U-shaped $(0.5,0.5)$, bell-shaped $(5,5)$,
left-skewed $(5,2)$ and right-skewed $(2,5)$; see Figure~\ref{fig:beta} for a graphical
illustration.

%

\begin{figure}[b]

\includegraphics{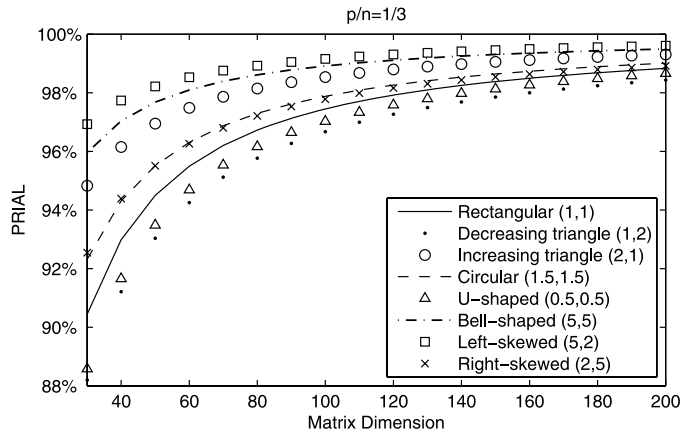}

\caption{Performance of the nonlinear shrinkage with beta densities. The various
curves correspond to different shapes of the population spectral
density. The support of the population spectral density is $[1,10]$.}\label{fig:shape}
\end{figure}

For every one of these eight beta densities, we take the population
eigenvalues to be equal to
\[
1+9F_{(\alpha,\beta)}^{-1} \biggl(\frac{i}{p}-\frac{1}{2p} \biggr),\qquad i=1,\ldots,p.
\]
The concentration ratio $\widehat{c}_n=p/n$ is equal to $1/3$. For
various values of $p$ between $30$ and $200$, we run $1000$
simulations with normally distributed variables. The PRIAL of the
nonlinear shrinkage estimator $\widehat{S}_n$
is plotted in Figure~\ref{fig:shape}.

%

As in all the other simulations presented above, the PRIAL of the
nonlinear shrinkage estimator always exceeds $88\%$, and more often
than not exceeds $95\%$. To preserve the clarity of the
picture, we do not report the PRIALs of the oracle and of the linear shrinkage
estimator; but as usual, the nonlinear shrinkage estimator ranked
between them.

\subsection{Fixed-dimension asymptotics}

\begin{figure}[b]

\includegraphics{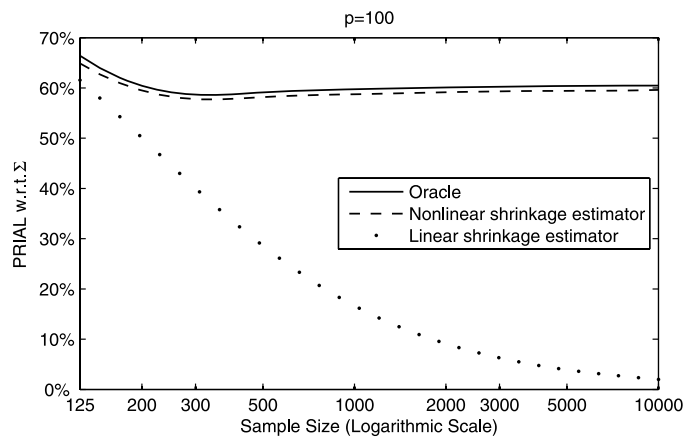}

\caption{Fixed-dimension asymptotics. $20\%$ of population eigenvalues are
equal to $1$, $40\%$ are equal to $3$ and $40\%$ are equal to
$10$. Variables are normally distributed. Every point is the result of
$1000$ Monte Carlo simulations.}\label{fig:fixeddim}
\end{figure}

Finally, we report a set of Monte Carlo simulations that departs from
the large-dimensional asymptotics assumption under which the nonlinear
shrinkage estimator $\widehat{S}_n$ was derived. The goal is to
compare it against the sample covariance matrix $S_n$ in the setting
where $S_n$ is known to have certain optimality properties (at least
in the normal case): traditional asymptotics, that is, when the number
of variables $p$ remains fixed while the sample size $n$ goes to
infinity. This gives as much advantage to the sample covariance matrix
as it can possibly have. We fix the dimension $p=100$ and let the
sample size $n$ vary from $n=125$ to $n=10\mbox{,}000$. In practice, very few
applied researchers are fortunate enough to have as many as $n=10\mbox{,}000$
i.i.d. observations, or a concentration ratio $c=p/n$ as low as
$0.01$. The respective PRIALs of $S_n^{or}$, $\widehat{S}_n$ and
$\overline{S}_n$ are plotted in Figure~\ref{fig:fixeddim}.

%

One crucial difference with all the previous simulations is that the
target for the PRIAL is no longer $S_n^*$, but instead the population
covariance matrix~$\Sigma$ itself, because now~$\Sigma$ can be
consistently estimated. Note that, since the matrix dimension is
fixed,~$\Sigma_n$ does not change with $n$; therefore, we can drop the
subscript $n$. Thus, in this subsection only, the definition of the
PRIAL is given by
\[
\mathrm{PRIAL} \equiv\mathrm{PRIAL}(\widehat\Sigma_n) \equiv
100\times
\biggl\{1-\frac{\E[ \Vert\widehat{\Sigma}_n-\Sigma\Vert^2 ]}
{\E[ \Vert S_n-\Sigma\Vert^2 ]} \biggr\} \%,
\]
where $\widehat\Sigma_n$ is an arbitrary estimator of
$\Sigma$. By definition, the PRIAL of $S_n$ is 0\%
while the PRIAL of $\Sigma$ is 100\%.

In this setting, \citet{ledoit:wolf:2004a} acknowledge that the
improvement of the linear shrinkage estimator over the sample
covariance matrix vanishes asymptotically, because the optimal linear
shrinkage intensity vanishes. Therefore it should be no surprise that
the PRIAL of $\overline{S}_n$ goes to zero in Figure~\ref{fig:fixeddim}. Perhaps more surprising is the continued ability
of the oracle and the nonlinear shrinkage estimator to improve by
approximately $60\%$ over the sample covariance matrix, even for a
sample size as large as $n=10\mbox{,}000$, and with no sign of abating as $n$
goes to infinity. This is an encouraging result, as our simulation
gave every possible advantage to the sample covariance matrix by
placing it in the asymptotic conditions where it possesses well-known
optimality properties, and where the earlier linear shrinkage
estimator of \citet{ledoit:wolf:2004a} is most disadvantaged.

\begin{figure}

\includegraphics{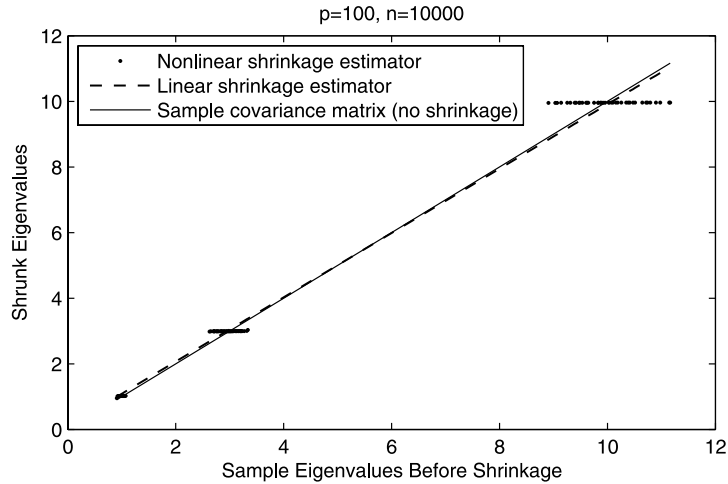}

\caption{Nonlinear shrinkage under fixed-dimension aymptotics. $20\%$ of
population eigenvalues are equal to $1$, $40\%$ are equal to $3$ and
$40\%$ are equal to $10$. $p=100$ and $n=10\mbox{,}000$. The oracle is not
shown because it is virtually identical to the nonlinear shrinkage estimator.}\label{fig:separation}
\end{figure}

Intuitively, this is because the oracle shrinkage formula becomes more
and more nonlinear as $n$ goes to infinity for fixed
$p$. \citet{bai:silverstein:1998} show that the sample covariance
matrix exhibits ``spectral separation'' when the concentration ratio
$p/n$ is sufficiently small. It means that the sample eigenvalues coalesce
into clusters, each cluster corresponding to a Dirac
of population eigenvalues. Within a given cluster, the
smallest sample eigenvalues need to be nudged upward, and the largest
ones downward, to the average of the cluster. In other words: full
shrinkage within clusters, and no shrinkage between clusters. This is
illustrated in Figure~\ref{fig:separation},
which plots a~typical simulation result for $n=10\mbox{,}000$.\setcounter{footnote}{1}\footnote{For
enhanced ability to distinguish linear shrinkage from the sample
covariance matrix, we plot the two uninterrupted lines, even though
the sample eigenvalues lie in three disjoint intervals (as can be
seen from nonlinear shrinkage).}

%

By detecting this intricate pattern automatically, that is, by
discovering where to shrink and where not to shrink, the nonlinear
shrinkage estimator $\widehat{S}_n$ showcases its ability to generate
substantial improvements over the sample covariance matrix even for
very low concentration ratios.


\subsection{Additional Monte Carlo simulations}
\label{ss:additional}

\subsubsection{Comparisons with other estimators}

So far, we have compared the nonlinear shrinkage
estimator $\widehat S_n$ only to the linear shrinkage estimator
$\overline
S_n$ and the oracle estimator $S_n^{or}$ to keep the
resulting figures concise and legible.

It is of additional interest to compare the nonlinear
shrinkage estimator also to some other estimators from the
literature. To this end we consider the following set of estimators:

\begin{itemize}
\item The estimator of \citet{stein:1975};
\item The estimator of \citet{haff:1980};
\item The estimator recently proposed by \citet{won:et:al:2009}. This
estimator is based
on a maximum likelihood approach, assuming normality, with an
explicit constraint on the condition number of the covariance
matrix. The\vadjust{\goodbreak} resulting estimator turns out to be a nonlinear
shrinkage estimator as well: all ``small'' sample eigenvalues are
brought up to a lower bound, all ``large'' sample eigenvalues are
brought down to an upper bound, and all ``intermediate'' sample
eigenvalues are left unchanged.

Therefore, the corresponding
transformation from sample eigenvalues to shrunk eigenvalues is
step-wise linear: first flat, then a 45-degree line, and then flat
again. The upper and lower bounds are determined by the desired
constraint on the condition number $\kappa$. If such an explicit constraint
is not available from a priori information, a suitable
constraint number $\widehat\kappa$ can be computed in a
data-dependent fashion by a
$K$-fold cross-validation method, which is the method we
use.\footnote{We are grateful to Joong-Ho Won for supplying us with
corresponding Matlab code.}

In particular, the cross-validation method selects $\widehat\kappa$
by optimizing over a finite grid $\{\kappa_1, \kappa_2, \ldots,
\kappa_L\}$ that has to be supplied by the user. To this end we
choose $L = 10$ and the $\kappa_l$ log-linearly spaced between 1 and
$\kappa(S_n)$, for $l = 1, \ldots, L$; here $\kappa(S_n)$ denotes
the condition number of the sample covariance matrix.
More precisely, for $l = 1, \ldots, L$, $\kappa_l \equiv\exp(\omega_l)$,
where
$\{\omega_1, \omega_2, \ldots, \omega_L\}$ is the equally-spaced
grid with
$\omega_1 \equiv0$ and $\omega_L \equiv\log(\kappa(S_n))$.
\item The cross-validation version of the nonlinear shrinkage
estimator $\widehat S_n$; see Remark~\ref{r:cv}.
\end{itemize}

\begin{figure}

\includegraphics{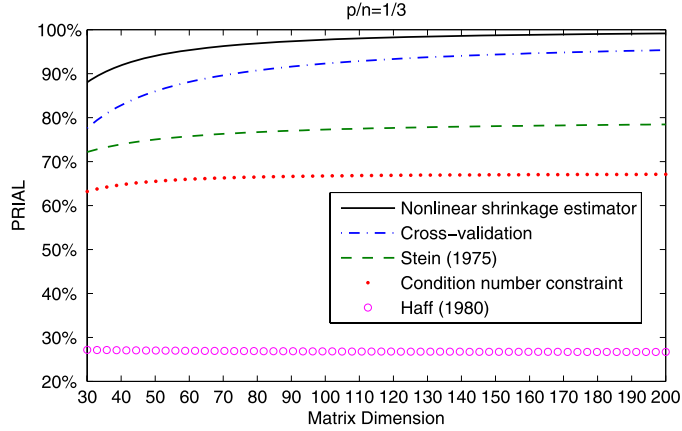}

\caption{Comparison of various estimators. $20\%$ of population eigenvalues are
equal to~$1$, $40\%$ are equal to $3$ and $40\%$ are equal to $10$.
Every point is the result of $1000$ Monte Carlo
simulations.}\label{fig:convergence:various}\vspace*{3pt}
\end{figure}

%
\begin{figure}[b]

\includegraphics{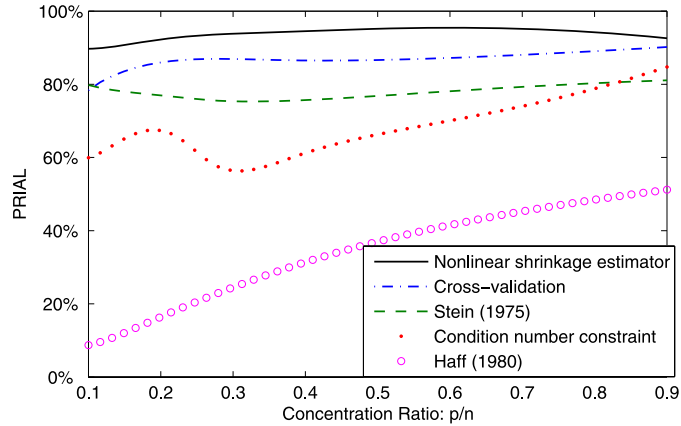}

\caption{Effect of varying the concentration ratio $\widehat{c}_n=p/n$. $20\%$
of population eigenvalues are equal to $1$, $40\%$ are equal to $3$ and
$40\%$ are equal to $10$. Every point is the
result of $1000$ Monte Carlo simulations.}\label{fig:concentration:various}
\end{figure}

We repeat the simulation exercises of
Sections~\ref{sub:MC_convergence}--\ref{sub:MC_dispersion}, replacing
the oracle estimator and the linear shrinkage estimator with the
above set of other estimators.
The respective PRIALs of the various estimators are plotted in
Figures~\ref{fig:convergence:various}--\ref{fig:dispersion:various}.

One can see that the nonlinear shrinkage estimator $\widehat S_n$
outperforms all other estimators, with the cross-validation version of
$\widehat S_n$ in second place, followed by the estimators of
\citet{stein:1975}, \citet{won:et:al:2009} and
\citet{haff:1980}.\looseness=1

\subsubsection{Comparisons based on a different loss function}

So far, the PRIAL has been based on the loss function
\[
L^{Fr}(\widehat\Sigma_n, \Sigma_n) \equiv
\|\widehat\Sigma_n - \Sigma_n \|^2.
\]
It is of additional interest to add some
comparisons based on a different loss function. To this end we use the
scale-invariant loss function proposed by \citet{james:stein:1961}, namely
%
\begin{equation} \label{e:loss-stein}
L^{JS} (\widehat\Sigma_n, \Sigma_n) \equiv \operatorname{trace} (\widehat
\Sigma_n \Sigma_n^{-1} ) - \log\det( \widehat
\Sigma_n \Sigma_n^{-1} ) - p.
\end{equation}

We repeat the simulation exercises of
Sections~\ref{sub:MC_convergence}--\ref{sub:MC_dispersion},
replacing $L^{Fr}$ with~$L^{JS}$. The respective PRIALs of $S_n^{or}$,
$\widehat{S}_n$, and
$\overline{S}_n$ are plotted in
Figures~\ref{fig:convergence:stein}--\ref{fig:dispersion:stein}.

%
\begin{figure}

\includegraphics{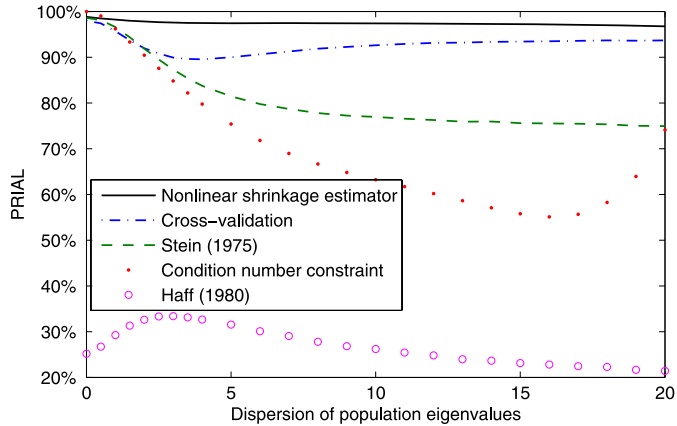}

\caption{Effect of varying the dispersion of population eigenvalues. $20\%$ of
population eigenvalues are equal to $1$, $40\%$ equal to $1+2 d/9$
and $40\%$ equal to $1+d$, where the dispersion parameter $d$ varies
from $0$ to $20$. $p=100$ and $n=300$.
Every point is the result of $1000$ Monte Carlo simulations.}\label{fig:dispersion:various}
\end{figure}

One can see that the results do not change much qualitatively. If
anything, the
comparisons are now even more favorable to the nonlinear shrinkage
estimator, in particular when comparing
Figure~\ref{fig:dispersion} to Figure~\ref{fig:dispersion:stein}.

%
\begin{figure}

\includegraphics{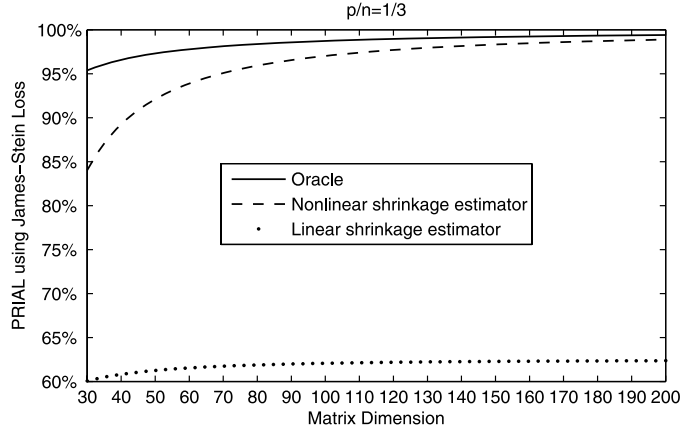}

\caption{Comparison of the nonlinear vs. linear shrinkage estimators. $20\%$ of
population eigenvalues are equal to $1$, $40\%$ are equal to $3$ and
$40\%$ are equal to $10$. The PRIALs are based on the
James--Stein (\protect\citeyear{james:stein:1961}) loss function \protect\eqref{e:loss-stein}. Every point is the
result of $1000$ Monte Carlo
simulations.}\label{fig:convergence:stein}\vspace*{3pt}
\end{figure}

%
\begin{figure}[b]
\vspace*{3pt}
\includegraphics{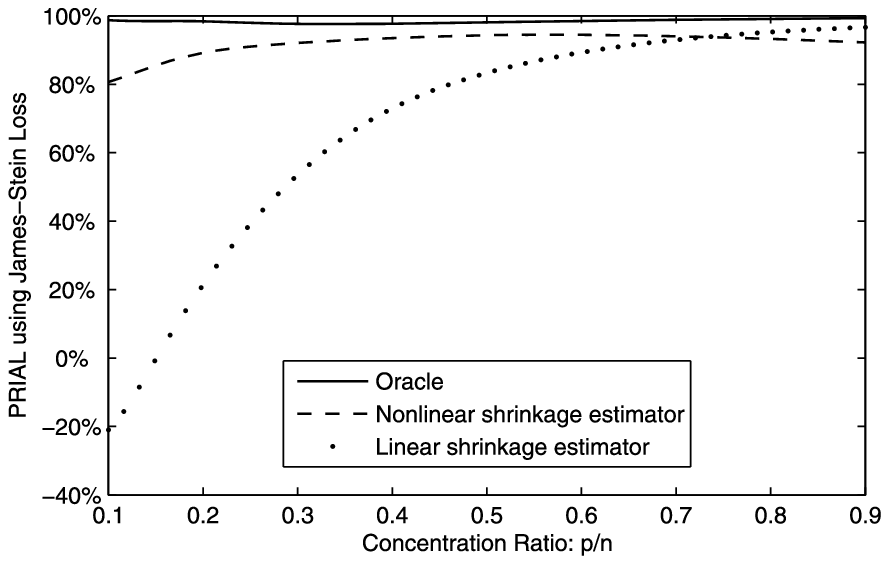}

\caption{Effect of varying the concentration ratio $\widehat{c}_n=p/n$. $20\%$
of population eigenvalues are equal to $1$, $40\%$ are equal to $3$ and
$40\%$ are equal to $10$. The PRIALs are based on the
James--Stein (\protect\citeyear{james:stein:1961}) loss function \protect\eqref{e:loss-stein}. Every point is the
result of $1000$ Monte Carlo simulations.}\label{fig:concentration:stein}
\end{figure}

%
\begin{figure}

\includegraphics{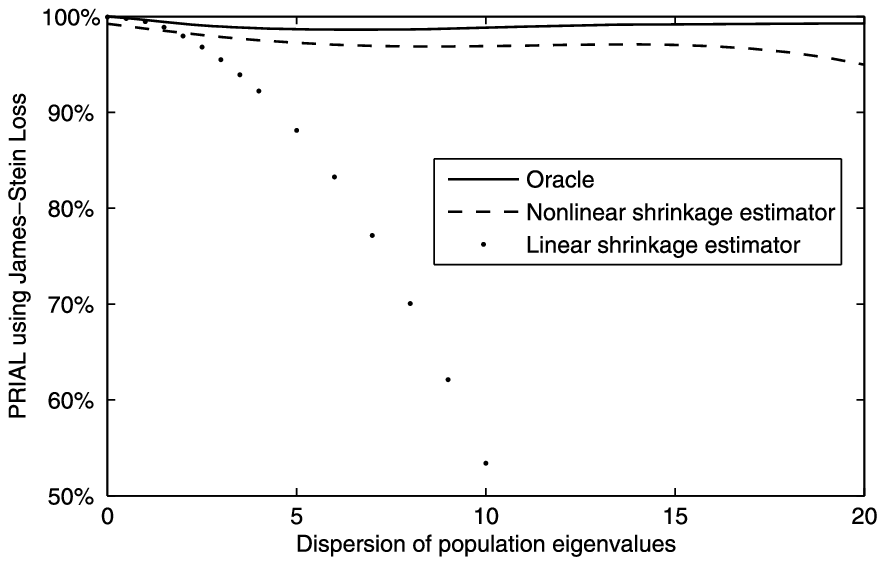}

\caption{Effect of varying the dispersion of population eigenvalues. $20\%$ of
population eigenvalues are equal to $1$, $40\%$ equal to $1+2 d/9$
and $40\%$ equal to $1+d$, where the dispersion parameter $d$ varies
from $0$ to $20$. $p=100$ and $n=300$. The PRIALs are based on the
James and Stein (\protect\citeyear{james:stein:1961}) loss function \protect\eqref{e:loss-stein}.
Every point is the result of $1000$ Monte Carlo
simulations.}\label{fig:dispersion:stein}\vspace*{6pt}
\end{figure}

\section{Conclusion}
\label{sec:conclusions}

Estimating a large-dimensional covariance matrix is
a~very important and challenging problem. In the absence of additional
information concerning the structure of the true covariance matrix, a
successful approach consists of appropriately shrinking the
sample eigenvalues, while retaining the sample eigenvectors.
In particular, such shrinkage estimators enjoy the desirable property
of being rotation-equivariant.

In this paper, we have
extended the linear approach of \citet{ledoit:wolf:2004a} by applying a
nonlinear transformation to the sample eigenvalues. The
specific
transformation suggested is motivated by the oracle estimator of
\citet{ledoit:peche:2011}, which in turn was derived by studying the
asymptotic behavior of the finite-sample optimal rotation-equivariant
estimator (i.e., the estimator with the rotation-equivariant property
that is closest to
the true covariance matrix when distance is measured by
the Frobenius norm).

The oracle estimator
involves the Stieltjes transform of the limiting spectral distribution of
the sample eigenvalues, evaluated at various points on the real
line. By finding a way to consistently estimate these quantities, in a
uniform sense, we have been able to construct a \textit{bona fide} nonlinear
shrinkage estimator that is asymptotically equivalent to the
oracle.

Extensive Monte Carlo studies have demonstrated the improved
finite-sample properties of our nonlinear shrinkage estimator compared
to the sample covariance matrix and the linear shrinkage estimator of
\citet{ledoit:wolf:2004a}, as well as its fast convergence to the
performance of the oracle.
In particular, when the sample size is very large
compared to the dimension, or the population eigenvalues are very
dispersed, the nonlinear shrinkage estimator still yields a significant
improvement over the sample covariance matrix, while the linear
shrinkage estimator no longer does.

Many statistical applications require an estimator of the inverse of the
covariance matrix, which is called the precision matrix. We have
modified our nonlinear shrinkage approach to this alternative problem,
thereby constructing a direct estimator of the precision matrix.
Monte Carlo studies have confirmed that this estimator yields a
sizable improvement over the indirect method of simply inverting the
nonlinear shrinkage estimator of the covariance matrix itself.

The scope of this paper is limited to the case where the matrix dimension
is smaller than the sample size. The other case, where the matrix
dimension exceeds
the sample size, requires certain modifications in the mathematical treatment,
and is left for future research.

\section*{Acknowledgments}
 We would like to thank two anonymous referees
for valuable comments, which have resulted in an improved exposition of
this paper.


\begin{supplement}
\stitle{Mathematical proofs}
\slink[doi]{10.1214/12-AOS989SUPP} 
\sdatatype{.pdf}
\sfilename{AOS989\_supp.pdf}
\sdescription{This supplement contains detailed proofs of
all mathematical results.}
\end{supplement}


%

\printaddresses

\end{document}